\documentclass{siamltex1213}
\usepackage{graphicx,a4wide}
\usepackage[intlimits]{amsmath}
\usepackage{amsxtra, amssymb, latexsym}%amsthm,
\usepackage[english]{babel}
\numberwithin{equation}{section}
\numberwithin{figure}{section}

\usepackage{tikz}
\usepackage{color}

\begin{document}

\newtheorem{claim}[theorem]{Claim}
\newtheorem{postulate}[theorem]{Postulate}
\newtheorem{assumption}[theorem]{Assumption}
\newtheorem{remark}[theorem]{Remark}
\newtheorem{example}[theorem]{Example}

% The following commands have been added by Sander Hille in order to avoid problems with the package bbm:
\newcommand{\mathbbm}[1]{{#1\!\!#1}}
\newcommand{\ind}{{\mathbbm{1}}}

% Further new commands:
\newcommand{\CM}{\mathcal{M}}
\newcommand{\CMc}{\overline{\CM}}
\newcommand{\smfrac}[2]{\mbox{$\frac{#1}{#2}$}}
\newcommand{\geqs}{\geqslant}
\newcommand{\leqs}{\leqslant}
\newcommand{\nn}{{(n)}}
\newcommand{\mm}{{(m)}}

% === end of commands added ===

\newcommand{\Lip}{{\mathrm{L}}}
\newcommand{\BL}{{\mathrm{BL}}}
\newcommand{\TV}{{\mathrm{TV}}}
\newcommand{\FM}{{\mathrm{FM}}}
\newcommand{\pair}[2]{\left\langle #1 , #2 \right\rangle}
\newcommand{\Int}[4]{\int_{#1}^{#2}\! #3 \, #4}
\newcommand{\map}[3]{#1 : #2 \rightarrow #3}
\def\R{\mathbb{R}}
\def\Rp{\mathbb{R}^+}
\def\N{\mathbb{N}}
\def\Np{\mathbb{N}^+}
\def\eps{\varepsilon}
\def\aeps{a}
\def\K{\mathcal{K}}

\title{Measure-valued mass evolution problems with flux boundary conditions and solution-dependent velocities}

\author{Joep H.M.~Evers\thanks{Department of Mathematics, Simon Fraser University, Burnaby, Canada, and Department of Mathematics and Statistics, Dalhousie University, Halifax, Canada. Corresponding author; email: \texttt{jevers@sfu.ca}.} \and  Sander C.~Hille\thanks{Mathematical Institute, Leiden University, P.O. Box 9512, 2300 RA, Leiden, The Netherlands.}% Corresponding author; email: \texttt{shille@math.leidenuniv.nl}}  
\and Adrian Muntean\thanks{Department of Mathematics and Computer Science, Karlstad University, Sweden.}}

%\author{Joep Evers, Sander Hille\thanks{Corresponding author}, Adrian Muntean}

% \subjclass[2000]{...}
% \keywords{measure-valued solution, flux boundary condition, boundary layer, singular limit, nonsmooth analysis}

%\date{\today}
\maketitle

\pagestyle{myheadings}
\thispagestyle{plain}
\markboth{J.H.M.~EVERS, S.C.~HILLE AND A.~MUNTEAN}{MEASURE-VALUED EQUATIONS WITH SOLUTION-DEPENDENT VELOCITY}

\begin{abstract}
In this paper we prove well-posedness for a measure-valued continuity equation with solution-dependent velocity and flux boundary conditions, posed on a bounded one-dimensional domain. We generalize the results of \cite{EvHiMu_JDE} to settings where the dynamics are driven by interactions. In a forward-Euler-like approach, we construct a time-discretized version of the original problem and employ the results of \cite{EvHiMu_JDE} as a building block within each subinterval. A limit solution is obtained as the mesh size of the time discretization goes to zero. Moreover, the limit is independent of the specific way of partitioning the time interval $[0,T]$.\\
This paper is partially based on results presented in \cite[Chapter 5]{Evers_PhD}, while a number of issues that were still open there, are now resolved.
\end{abstract}



%46E27  	(Functional Analysis) Spaces of measures
%28A33  	Spaces of measures, convergence of measures
%35Q91 PDEs in connection with game theory, economics, social and behavioral sciences
% 34A12 Initial value problems, existence, uniqueness, continuous dependence and continuation of solutions
%35L65  Conservation laws

%\tableofcontents

\section{Introduction}
A considerable amount of recent mathematical literature has been devoted to evolution equations formulated in terms of measures. Such equations are used to describe systems that occur in e.g.~biology (animal aggregations \cite{CarrilloFornasierRosadoToscani,CanCarRos}, crowds of pedestrians \cite{CrisPiccTosinBook}, structured populations \cite{DG, Gwiazda1,Gwiazda2,Ackleh1}) and material science (defects in metallic crystals \cite{vMeurs14}). Many interesting and relevant scenarios take place in \textit{bounded} domains. Apart from the examples mentioned above, these include intracellular transport processes, cf.~\cite[Section 1]{EversMBE}, and also manufacturing chains \cite{Schleper}. However, most works that deal with well-posedness of measure-valued equations and properties of their solutions treat these equations in the full space, see for instance also \cite{BenzoniColomboGwiazda, CarrDiFranFigLauSlep11, TosinFrasca,CrippaLecureux, Colombo_Crippa_Spirito}. The present work explicitly focuses on bounded domains and the challenge of defining mathematically and physically `correct' boundary conditions.\\
\\
In \cite{EvHiMu_JDE}, we derived boundary conditions for a one-dimensional measure-valued transport equation on the unit interval $[0,1]$ with prescribed velocity field $v$. A short-hand notation for this equation is:
\begin{equation}\label{eq:main equation}
\frac{\partial}{\partial t}\mu_t + \frac{\partial}{\partial x} (v\,\mu_t) = f\cdot \mu_t.
\end{equation}
We focused on the well-posedness of this equation, in the sense of \emph{mild solutions}, and the convergence of solutions corresponding to a sequence $(f_n)_{n\in\N}$ in the right-hand side. Some specific choices for $(f_n)_{n\in\N}$ represent for instance effects in a boundary layer that approximate, as $n\to\infty$, sink or source effects localized on the boundary. The boundary layer corresponds to the regions in $[0,1]$ where the functions $f_n$ are nonzero.\\
\\
%can, for instance, be viewed as a way to approximate sink or source effects localized on the boundary by effects in a boundary layer. The regions in $[0,1]$ where functions $f_n$ (and their limit as $n\to\infty$) are nonzero plays an important role in this respect.\\
%\\
There are several reasons \textit{why} we consider mild solutions rather than weak solutions. First of all, the mild formulation in terms of the variation of constants formula -- see \eqref{eq:VOC} -- follows directly from a probabilistic interpretation, as was shown in \cite[Section 6]{EvHiMu_JDE}. Therefore the choice for mild solutions is justified by a modelling argument. Secondly, usually uniqueness of weak solutions cannot be expected to hold, while mild solutions are unique when the perturbation ($\mu\mapsto f\cdot\mu$) is Lipschitz. In \cite{EvHiMu_JDE}, where the perturbation even has discontinuities, we still obtain uniqueness of the mild solution. This is one of the main results of \cite{EvHiMu_JDE}. In the works \cite{Ackleh1, Gwiazda1, Gwiazda2,Gwiazda-Jamroz_ea:2012} %<---removed Ackleh-Fitzpatrick-Thieme, paper not available
a specific weak solution is constructed that is precisely the mild solution that we obtain by different means. %uniqueness is obtained by constructing the mild solution.
Finally, there is a technical advantage of using mild solutions. Most of our estimates are in terms of the dual bounded Lipschitz norm $\|\cdot\|_\BL^*$, that will be introduced in \S \ref{sec: basics meas th}. Because test functions do not appear explicitly, our calculations are often simpler than when weak solutions are considered. Moreover, our estimates are in fact uniform over test functions in a bounded set.\\
\\
In the present work, we propose and investigate a procedure to generalize the former results to include velocity fields that depend on the solution itself. Such generalization makes it possible to model in a bounded domain the dynamics governed by interactions between the `particles'; in particular we will be concerned with interaction terms of convolution type that are given by a weighted average over the whole population.\\
\\
The results in this paper hold for a source-sink right-hand side that is based on a function $f$ that is an element of the space $\BL([0,1])$ of bounded Lipschitz functions on $[0,1]$. In \cite{EvHiMu_JDE}, we worked with $\map{f}{[0,1]}{\R}$ that is \textit{piecewise} bounded Lipschitz, though. Hence, here we are able to describe absorption in a boundary layer, but not yet absorption on the boundary alone. In the discussion section of this paper, see \S \ref{sec: disc f}, we comment on the possibilities to extend our results to $f$ that is piecewise bounded Lipschitz.\\
\\
%Mention the approach of \cite{Gwiazda1}, Prop 3.11, Ron mentions this (no uniqueness)
%cite related approaches: Piccoli-Rossi 2011 ``Transport equation with nonlocal velocity...'', Ca{\~n}izo-Carrillo-Cuadrado?,  Colombo-Guerra DCDS 2009 Euler polygonals in general metric spaces (cf. Lorentz workshop),
We consider \eqref{eq:main equation} for velocity fields that are no longer \textit{fixed} elements of $\BL([0,1])$. Instead of $v$, we write $v[\mu]$ for the velocity field that depends functionally on the measure $\mu$. The transport equation on $[0,1]$ becomes
\begin{equation}\label{eqn:main equation nonlin}
\frac{\partial}{\partial t}\mu_t + \frac{\partial}{\partial x} (v[\mu_t]\,\mu_t ) = f\cdot \mu_t.
\end{equation}
%,The fact that we \textit{extend} the former results is for instance reflected by the fact that we will assume that $v[\mu]$ is still an element of $\BL([0,1])$, in this case for each $\mu\in\CM^+([0,1])$. Further requirements of the velocity field are given in Assumption \ref{ass: v properties}.\\
The aim of this paper is to ensure the well-posedness of \eqref{eqn:main equation nonlin}, in a suitable sense. Because \eqref{eqn:main equation nonlin} is a nonlinear equation, establishing well-posedness is not straightforward. Here, we employ a \textit{forward-Euler-like approach} that builds on the fundamentals constructed in \cite{EvHiMu_JDE}. We partition the time interval $[0,T]$ and fix the velocity on each subinterval. That is, restricted to a subinterval, the velocity depends only on the spatial variable and not on the solution measure. Within each subinterval the measure-valued solution evolves according to the fixed velocity and the evolution fits in the framework set in \cite{EvHiMu_JDE}. A more detailed description of our approach is given in \S \ref{sec: gen meas-dep v}. We decrease the mesh size in the partition of $[0,T]$ and estimate the difference between Euler approximations. The main result of this paper is the fact that this procedure converges.\\
\\
%We were inspired by \cite{ColomboGuerra} and \cite{Hoogwater-thesis}, the latter of which was in turn inspired by \cite{Gwiazda1}.
A forward-Euler scheme similar to ours is used in \cite{PiccoliRossi} for measures absolutely continuous with respect to the Lebesgue measure. Their results are extended to general measures in \cite[Chapter 7]{CrisPiccTosinBook}. The difference between their work and ours is twofold: they use the Wasserstein distance and they work in unbounded domains.\\
The references that directly inspired us are \cite{ColomboGuerra, Hoogwater-thesis,Gwiazda1}. The approach presented in this paper deviates from \cite{Hoogwater-thesis}, since we restrict ourselves to evolution on the interval $[0,1]$, while \cite{Hoogwater-thesis} considers $[0,\infty)$. Furthermore, our regularity conditions on the velocity -- given in Assumption \ref{ass: v properties} -- are weaker than in \cite{Hoogwater-thesis}; cf.~Remark \ref{rem:estimates on ball i.o. unif Hoogwater}. Moreover, \cite{Hoogwater-thesis} restricts to velocity fields that point inwards at $0$. In this way, no mass is allowed to flow out of the domain $[0,\infty)$. In our approach, the fact that the flow is stopped at the boundary is encoded in the semigroup $(P_t)_{t\geqs0}$, irrespective of the sign of the velocity there; cf.~\S \ref{sec: prop stopped flow}. We consider it too restrictive to have a condition on the sign of the velocity at $0$ or $1$; in practice it is very difficult to make sure that such condition is satisfied when the velocity $v[\mu]$ depends on the solution (like in e.g.~Example \ref{ex: meas-dep v}).\\
\\
In this paper we limit our attention to a one-dimensional state space, $[0,1]$, because in this case the (global) Lipschitz continuous dependence of the stopped flow on the  time-invariant velocity field $v$ is a rather straightforward property (see Section \ref{sec: prop stopped flow}, Lemma \ref{lem: Pv mu - Pv' mu}). In higher dimensional (bounded) state spaces this seems much more delicate to establish. We comment on this in more detail after the proof of Lemma \ref{lem: Pv mu - Pv' mu}. One should note however, that the results on convergence of the forward-Euler-like approach that we present do not depend on the dimensionality other than through the mentioned Lipschitzian property as presented in Lemma \ref{lem: Pv mu - Pv' mu}.
\\
%As said, on each of the subintervals in the forward-Euler-like approach, we build on the results of Chapter \ref{ch:bc pres v}. For convenience we therefore restate the relevant results. Let $v:[0,1]\to\R$ be a given velocity field that is bounded Lipschitz, and let $P_t$ be the semigroup on $\CM([0,1])$ that corresponds to the stopped flow induced by $v$. The function $f$ is a {\em piecewise bounded Lipschitz function} on $[0,1]$ with finitely many discontinuities. We associate to $f$ the perturbation map $F_f:\CM([0,1])\to\CM([0,1]): \mu\mapsto f\cdot \mu$.\\
%\\
%In Chapter \ref{ch:bc pres v}, we show that for each given initial measure $\nu\in\CM^+([0,1])$, the variation of constants formula\index{variation of constants formula}
%\begin{equation}\label{eq:VOC summ}
%\mu_t = P_t\,\nu \ +\ \int_0^t P_{t-s}F_f(\mu_s)\, ds\qquad \mbox{for all}\ t\in [0,T].
%\end{equation}
%has a unique solution $\mu_\cdot:[0,T]\to\CM^+([0,1])_\BL$ that is $\|\cdot\|_\TV$-bounded; see Corollary \ref{cor:exist uniq}. Moreover, we proved in Proposition \ref{prop:cont dependence} that the solution depends continuously on initial data.

This paper is organized as follows. Within each subinterval of the Euler approximation the dynamics are given by a fixed velocity, like in \cite{EvHiMu_JDE}. Therefore, we start in \S \ref{sec: summ tech} by collecting the results of \cite{EvHiMu_JDE} that we require here: a number of properties of the semigroup $(P_t)_{t\geqs0}$ and of the solution operator, called $(Q_t)_{t\geqs0}$. %In particular, we provide a Lipschitz estimate for semigroups that follow from two different velocity fields; see Lemma \ref{lem: Pv mu - Pv' mu}.
The forward-Euler-like approach to construct solutions is introduced in \S \ref{sec: gen meas-dep v}, where we also state the main results of this paper: Theorems \ref{thm: exist uniq nonlin} and \ref{thm: cont dep}, and Corollary \ref{cor: global exist uniq}. In plain words and combined into one pseudo-theorem, these results read:\\
\\
\textsc{Theorem.} \emph{The proposed forward-Euler-like approach converges as the mesh size of the time discretization goes to zero. The limit is independent of the specific way in which the time domain is partitioned. This approximation procedure yields existence and uniqueness of mild solutions to the nonlinear problem, and solutions depend continuously on initial data.
}\\
\\
A more precise formulation follows later. We prove these results in \S \ref{sec: proofs ch bc nonlin} using estimates between two Euler approximations of \eqref{eq:main equation}. In \S \ref{sec: discussion bc meas-dep v} we reflect on the achievements of this paper, discuss open issues and provide directions for further research.

\section{Preliminaries}\label{sec: summ tech}
This section contains a summary of the results obtained in \cite{EvHiMu_JDE} on which we shall build. Moreover, we mention the technical preliminaries needed for the arguments in this paper. 

\subsection{Basics of measure theory}\label{sec: basics meas th}
If $S$ is a topological space, we denote by $\CM(S)$ the space of finite Borel measures on $S$ and by $\CM^+(S)$ the convex cone of positive measures included in it.
%This cone defines a partial ordering on $\CM(S)$: $\mu\leqs \nu$ iff $\nu-\mu\in\CM^+(S)$. Clearly, $\mu\leqs \nu$ if and only if $\mu(E)\leqs \nu(E)$ for all Borel measurable $E\subset S$.
For $x\in S$, $\delta_x$ denotes the Dirac measure at $x$. Let
\begin{equation}\label{pairing}
\langle\mu,\phi\rangle :=\int_S \phi \,d\mu
\end{equation}
denote the natural pairing between measures $\mu\in\CM(S)$ and bounded measurable functions $\phi$. The {\em push-forward} or {\em image measure} of $\mu$ under Borel measurable $\Phi:S\to S$ is the measure $\Phi\#\mu$ defined on Borel sets $E\subset S$ by
\begin{equation}
(\Phi\# \mu)(E) := \mu\bigl(\Phi^{-1}(E)\bigr).\label{eqn: def push-forw chapter bc}
\end{equation}
One easily verifies that $\langle\Phi\# \mu,\phi\rangle=\langle\mu,\phi\circ\Phi\rangle$.\\
\\
We denote by $C_b(S)$ the Banach space of real-valued bounded continuous functions on $S$ equipped with the supremum norm $\|\cdot\|_\infty$. The {\em total variation norm} $\|\cdot\|_\TV$ on $\CM(S)$ is defined by
\begin{equation*}
\|\mu\|_{\TV}:= \sup\left\{\pair{\mu}{\phi}\,\Big|\, \phi\in C_b(S),\ \|\phi\|_\infty\leqslant 1 \right\}.
\end{equation*}
It follows immediately that for $\Phi:S\to S$ continuous, $\|\Phi\#\mu\|_\TV\leqs\|\mu\|_\TV$.
%
%The total variation norm is too strong for our application, since $\|\delta_x-\delta_y\|_\TV = 2$ if $x\neq y$. The natural topology to consider is the weak topology induced by $C_b(S)$ through the pairing \eqref{pairing}. In this topology $x\mapsto\delta_x:S\to\CM^+(S)$ is continuous.\\
%\\
In our setting, %$S=[0,1]$
$S$ is a Polish space (separable, completely metrizable topological space; cf.~\cite[p.~344]{Dudley}). It is well-established (cf. \cite{Dudley1,Dudley2}) that in this case the weak topology on $\CM(S)$ induced by $C_b(S)$ {\em when restricted to the positive cone} $\CM^+(S)$ is metrizable by a metric derived from a norm, e.g.~the Fortet-Mourier norm or the Dudley norm. The latter is also called the {\em dual bounded Lipschitz norm}, that we shall introduce now. To that end, let $d$ be a metric on $S$ that metrizes the topology, such that $(S,d)$ is separable and complete. Let $\BL(S,d)=\BL(S)$ be the vector space of real-valued bounded Lipschitz functions on $(S,d)$. For $\phi\in\BL(S)$, let
\begin{equation*}
|\phi|_\Lip  := \sup\left\{ \frac{|\phi(x)-\phi(y)|}{d(x,y)}\;\Big|\; x,y\in S,\ x\neq y\right\}
\end{equation*}
be its Lipschitz constant. Now
\begin{equation}\label{def:BL-norm}
\|\phi\|_\BL := \|\phi\|_\infty + |\phi|_\Lip
\end{equation}
defines a norm on $\BL(S)$ for which this space is a Banach space \cite{Fortet-Mourier:1953,Dudley1}. In fact, with this norm $\BL(S)$ is a Banach algebra for pointwise product of functions:
\begin{equation}\label{eqn:Banach algebra}
\|\phi\cdot\psi\|_\BL \le \|\phi\|_\BL\,\|\psi\|_\BL.
\end{equation}
Alternatively, one may define on $\BL(S)$ the equivalent norm%\footnote{See Lemma \ref{lemma: equivalence of norms} in Appendix \ref{app:equiv norms}.}
\begin{equation*}
\|\phi\|_\FM := \max\bigl(\|\phi\|_\infty\,,|\phi|_\Lip \bigr),
\end{equation*}
where `$\FM$' stands for `Fortet-Mourier' (see below). Let $\|\cdot\|_\BL^*$ be the dual norm of $\|\cdot\|_\BL$ on the dual space $\BL(S)^*$, i.e.~for any $x^*\in\BL(S)^*$ its norm is given by
\begin{equation*}
\|x^*\|_{\BL}^*:= \sup\left\{ |\pair{x^*}{\phi}|\ |\; \phi\in\BL(S),\ \|\phi\|_\BL\leqslant1\right\}.
\end{equation*}
The map $\mu\mapsto I_\mu$ with $I_\mu(\phi):=\pair{\mu}{\phi}$ defines a linear {\em embedding} of $\CM(S)$ into $\BL(S)^*$; see \cite[Lemma~6]{Dudley1}. Thus $\|\cdot\|^*_\BL$ induces a norm on $\CM(S)$, which is denoted by the same symbols. It is called the dual bounded Lipschitz norm or Dudley norm. Generally, $\|\mu\|_\BL^*\leqslant\|\mu\|_\TV$ for all $\mu\in\CM(S)$. For positive measures the two norms coincide:
\begin{equation}\label{TV norm is dual BL norm for pos measures}
\|\mu\|_{\BL}^*=\mu(S)=\|\mu\|_{\TV} \hspace{1 cm}\text{for all }\mu\in\CM^+(S).
\end{equation}
One may also consider the restriction to $\CM(S)$ of the dual norm $\|\cdot\|_\FM^*$ of $\|\cdot\|_\FM$ on $\BL(S)^*$. This yields an equivalent norm %\footnote{See Lemma \ref{lemma: equivalence of dual norms} in Appendix \ref{app:equiv norms}.}
on $\CM(S)$ that is called the Fortet-Mourier norm (see e.g.~\cite{Lasota-Myjak-Szarek,Zaharopol}):
\begin{equation}\label{eq:equivalence Dudely Fortet-Mourier}
\|\mu\|^*_\BL\leqs\|\mu\|^*_\FM\leqs 2\|\mu\|^*_\BL.
\end{equation}
This norm also satisfies $\|\mu\|^*_\FM\leqs\|\mu\|_\TV$, so \eqref{TV norm is dual BL norm for pos measures} holds for $\|\cdot\|_\FM^*$ too. Moreover (cf.~\cite[Lemma~3.5]{Hille-Worm:SF}), for any $x,y\in S$,
\begin{equation}\label{eq:Bl-norm difference Diracs}
\|\delta_x-\delta_y\|_\BL^* = \frac{2d(x,y)}{2+d(x,y)} \leqs \min(2,d(x,y)) = \|\delta_x-\delta_y\|_\FM^*.
\end{equation}
In general, the space $\CM(S)$ is not complete for $\|\cdot\|_\BL^*$. We denote by $\CMc(S)_\BL$ its completion, viewed as closure of $\CM(S)$ within $\BL(S)^*$. The space $\CM^+(S)$ is complete for $\|\cdot\|^*_\BL$, hence closed in $\CM(S)$ and $\CMc(S)_\BL$.\\
\\
The $\|\cdot\|_\BL^*$-norm is convenient also for integration. In Appendix~C of \cite{EvHiMu_JDE} some technical results about integration of measure-valued maps were collected. These will also be used in this paper.
The continuity of the map $x\mapsto\delta_x:S\to\CM^+(S)_\BL$ together with (C.2) in \cite{EvHiMu_JDE} yields the identity
\begin{equation}\label{eq:superposition of Diracs}
\mu = \int_S\delta_x\,d\mu(x)
\end{equation}
as Bochner integral in $\CMc(S)_\BL$; for basic results on Bochner integration, the reader is referred to e.g.~\cite{Diestel-Uhl}. The observation \eqref{eq:superposition of Diracs} will essentially link `continuum' (`$\mu$') and particle description (`$\delta_x$') for our equation on $[0,1]$.

\subsection{Properties of the stopped flow}\label{sec: prop stopped flow}
Let $v\in\BL([0,1])$ be fixed. We assume that a single particle (`individual') is moving in the domain $[0,1]$ deterministically, described by the differential equation for its position $x(t)$ at time $t$:
\begin{equation}\label{eq:indiv flow}
\left\{
  \begin{array}{l}
    \dot{x}(t)=v(x(t)), \\
    x(0)=x_0.
  \end{array}
\right.
\end{equation}
A solution to \eqref{eq:indiv flow} is unique, it exists for time up to reaching the boundary $0$ or $1$ and depends continuously on initial conditions. Let $x(\,\cdot\,;x_0)$ be this solution and $I_{x_0}$ be its maximal interval of existence. Define
\[
\tau_\partial(x_0) := \sup I_{x_0} \in [0,\infty],
\]
i.e.~$\tau_\partial(x_0)$ is the time at which the solution starting at $x_0$ reaches the boundary (if it happens) when $x_0$ is an interior point. Note that $\tau_\partial(x_0)=0$ when $x_0$ is a boundary point where $v$ points outwards, while $\tau_\partial(x_0)>0$ when $x_0$ is a boundary point where $v$ vanishes or points inwards.\\
\\
The {\em individualistic stopped flow} on $[0,1]$ associated to $v$ is the family of maps $\Phi_t:[0,1]\to[0,1]$, $t\geqs 0$,  defined by
\begin{equation}\label{individualistic flow Phi}
\Phi_t(x_0) := \begin{cases} x(t;x_0),& \quad \mbox{if}\ t\in I_{x_0},\\
x(\tau_\partial(x_0);x_0), & \quad \mbox{otherwise}.
\end{cases}
\end{equation}

%%%%%%%%%%%%%%%
To lift the dynamics to the space of measures, we define $\map{P_t}{\CM([0,1])}{\CM([0,1])}$ by means of the push-forward under $\Phi_t$: for all $\mu\in\mathcal{M}([0,1])$,
\begin{equation}\label{Def Pt push forward}
P_t\mu := \Phi_t \# \mu = \mu\circ\Phi_t^{-1};
\end{equation}
see \eqref{eqn: def push-forw chapter bc}. Clearly, $P_t$ maps positive measures to positive measures and $P_t$ is mass preserving on positive measures. Since the family of maps $(\Phi_t)_{t\geqs 0}$ forms a semigroup, so do the maps $P_t$ in the space $\CM([0,1])$. That is, $(P_t)_{t\geqs0}$ is a {\em Markov semigroup} on $\CM[0,1]$ (cf.~\cite{Lasota-Myjak-Szarek}). The basic estimate
\begin{equation}\label{eqn:TV norm Pt}
\|P_t\mu\|_\TV \leqs \|\mu\|_\TV
\end{equation}
holds for $\mu\in\CM([0,1])$.\\
\\
In the rest of this section we summarize those properties of $(P_t)_{t\geqs0}$ that are needed in this paper. We first recall Lemma~2.2 from \cite{EvHiMu_JDE}:\\

\begin{lemma}[See {\cite[Lemma 2.2]{EvHiMu_JDE}}]\label{lemma:lift is Lipschitz in time}
Let $\mu\in\CM([0,1])$ and $t,s\in\Rp$. Then
\begin{enumerate}
\item[(i)] $\| P_t\mu- P_s\mu\|_\BL^* \leqs \|v\|_\infty\,\|\mu\|_\TV \,|t-s|$.\label{lem: part BL bound t s}

\item[(ii)] $\| P_t\mu\|_\BL^* \leqs \max(1,|\Phi_t|_\Lip )\, \|\mu\|_\BL^* \leqs e^{|v|_\Lip t}\|\mu\|_\BL^*$.\label{lem: part BL bound Pmu}
\end{enumerate}
\end{lemma} %\begin{corollary}\label{cor: Pmu - Pnu}
%For all $\mu,\nu\in\CM([0,1])$, the following estimate holds:
%\begin{equation}
%\| P_t\mu-P_t\nu\|_\BL^* \leqs e^{|v|_\Lip t}\|\mu-\nu\|_\BL^* \,\,\,\text{for all }t\geqs0.
%\end{equation}
%\end{corollary}
%\begin{proof}
%Apply Part \ref{lem: part BL bound Pmu} of Lemma \ref{lemma:lift is Lipschitz in time} to the measure $\mu-\nu\in\CM([0,1])$.
%\end{proof}
\vspace{\baselineskip}
To distinguish between the semigroups on $\CM([0,1])$ associated to $v,v'\in\BL([0,1])$, we write $P^v$ and $P^{v'}$, respectively. Analogously, we distinguish between the semigroups $(\Phi^v_t)_{t\geqslant0}$ and $(\Phi^{v'}_t)_{t\geqslant0}$ on $[0,1]$ and between the intervals of existence $I^v_{x_0}$ and $I^{v'}_{x_0}$ associated to \eqref{eq:indiv flow}.\\

\begin{lemma}\label{lem: Pv mu - Pv' mu}
For all $\mu\in\CM([0,1])$, $v,v'\in\BL([0,1])$ and $t\in\Rp_0$
\begin{equation*}
\|P^{v}_t\mu - P^{v'}_t\mu\|^*_\BL \leqslant \|v-v'\|_\infty\,t\,\|\mu\|_\TV\,e^{L\,t},
\end{equation*}
where $L:=\min(|v|_\Lip\,,|v'|_\Lip)$.
\end{lemma}\\

\begin{proof}
For any $\phi\in\BL([0,1])$, we have
\begin{equation}\label{eqn: pair Pvt-Pv't phi Lipschitz}
|\pair{\phi}{P^v_t\mu-P^{v'}_t\mu}| = |\pair{\phi\circ \Phi^v_t - \phi\circ \Phi^{v'}_t}{\mu}| \leqslant |\phi|_\Lip\, \| \Phi^v_t - \Phi^{v'}_t  \|_\infty\,\|\mu\|_\TV,
\end{equation}
hence
\begin{equation}\label{eqn: BL norm Pv - Pv'}
\|P^v_t\mu-P^{v'}_t\mu\|^*_\BL \leqslant \| \Phi^v_t - \Phi^{v'}_t  \|_\infty\,\|\mu\|_\TV.
\end{equation}
Let $x\in[0,1]$.

\paragraph{Case 1: $t\in I^v_x\cap I^{v'}_x$}
\begin{align*}
\nonumber | \Phi^v_t(x) - \Phi^{v'}_t(x) | &= \left|\Int{0}{t}{v(\Phi^v_s(x)) - v'(\Phi^{v'}_s(x))}{ds}  \right| \\
&\leqslant |v|_\Lip \Int{0}{t}{|\Phi^v_s(x) - \Phi^{v'}_s(x)|}{ds} + \|v-v'\|_\infty\,t.
\end{align*}
Gronwall's Lemma yields
\begin{equation}
| \Phi^v_t(x) - \Phi^{v'}_t(x) | \leqslant \|v-v'\|_\infty\,t\, e^{|v|_\Lip\,t},\label{eqn: Pv - Pv' after Gronwall}
\end{equation}
for all $x\in[0,1]$. Due to the symmetry of \eqref{eqn: Pv - Pv' after Gronwall} in $v$ and $v'$, the same estimate \eqref{eqn: Pv - Pv' after Gronwall} can be obtained with $|v'|_\Lip$ instead of $|v|_\Lip$, and hence, we can write $\min(|v|_\Lip\,,|v'|_\Lip)$ in the exponent. This observation yields, together with \eqref{eqn: BL norm Pv - Pv'}, the statement of the lemma.

\paragraph{Case 2: $t\not\in I^v_x$}
We extend $v:[0,1]\to\R$ to $\bar{v}:\R\to\R$ by defining $\bar{v}(x):=v(0)$ if $x<0$ and $\bar{v}(x):=v(1)$ if $x>1$. Then $\bar{v}$ is a bounded Lipschitz extension of $v$ such that $\|\bar{v}\|_\infty=\|v\|_\infty$ and $|\bar{v}|_\Lip=|v|_\Lip$. Let $\Phi^{\bar{v}}_t:\R\to\R$ be the solution semigroup associated to the unique (global) solution to \eqref{eq:indiv flow} with $v$ replaced by $\bar{v}$ and with initial condition to be taken from the whole of $\R$. We extend $v'$ analogously to $\bar{v}'$.\\
Irrespective of whether $t\in I^{v'}_x$ or $t\not\in I^{v'}_x$, and whether in the latter case $\Phi^v_t(x)=\Phi^{v'}_t(x)$ or $\Phi^v_t(x)\neq\Phi^{v'}_t(x)$, the following estimate holds
\begin{equation}\label{eqn: est Pv-Pv'against extended}
| \Phi^v_t(x) - \Phi^{v'}_t(x) | \leqslant | \Phi^{\bar{v}}_t(x) - \Phi^{\bar{v}'}_t(x) |
\end{equation}
for all $x\in[0,1]$. estimate $| \Phi^{\bar{v}}_t(x) - \Phi^{\bar{v}'}_t(x) |$ using the same ideas as in \eqref{eqn: BL norm Pv - Pv'} and \eqref{eqn: Pv - Pv' after Gronwall} and obtain
\begin{equation*}
\|P^{\bar{v}}_t\mu - P^{\bar{v}'}_t\mu\|^*_\BL \leqslant \|\bar{v}-\bar{v}'\|_\infty\,t\,\|\mu\|_\TV\,\exp(\min(|\bar{v}|_\Lip\,,|\bar{v}'|_\Lip)\,t).
\end{equation*}
The statement of the lemma follows from the equalities $|\bar{v}|_\Lip=|v|_\Lip$, $|\bar{v}'|_\Lip=|v'|_\Lip$,  $\|\bar{v}-\bar{v}'\|_\infty= \|v-v'\|_\infty$ and Equation \eqref{eqn: est Pv-Pv'against extended}. The case $t\not\in I^{v'}_x$ is analogous.
\end{proof}
\\
\begin{remark}
The definition of stopped flow in state spaces of dimension two and higher and establishing elementary properties of its lift to measures is more delicate than the one-dimensional case presented above. Consider an open domain $\Omega\subset\R^n$, $n\geq 2$ (with sufficiently smooth boundary). Let $\overline{\Omega}$ be its closure and let $v\in\BL(\overline{\Omega}, \R^n)$ be a velocity field on $\overline{\Omega}$. Solutions to the initial value problem \eqref{eq:indiv flow} with $x_0\in\Omega$ still exist for some positive time, but in this higher dimensional setting it may happen that trajectories of the flow in $\overline{\Omega}$ defined by $v$ are partially contained in the boundary $\partial\Omega$ or even only `touch' $\partial\Omega$. So reaching the boundary in finite time is not equivalent to `leaving the domain'.

By means of the Metric Tietze Extension Theorem one can extend $v$ to $\bar{v}\in\BL(\R^n,\R^n)$ with preservation of Lipschitz constant and supremum norm, e.g. component-wise. Let $x_{\bar{v}}(t;x_0)$ be the corresponding solution starting at $x_0$ at $t=0$. The {\em exit time} or {\em stopping time} of the solution starting at $x_0$ could then be defined as
\[
\tau^v_\partial(x_0) := \inf\{t>0\,|\, x_{\bar{v}}(t;x_0)\in \R^n\setminus\overline{\Omega}\},
\]
with the convention that the infimum of the empty set is set to $+\infty$. One must show that this value is independent of the particular extension $\bar{v}$ that was chosen. $\tau^v_\partial(x_0)$ now replaces the similarly denoted termination time of the solution that was used above for the one-dimensional case. The stopped flow $\Phi^v_t$ can then be defined as in \eqref{individualistic flow Phi}. %Note that a naive definition `stop once on the boundary' will not yield a semigroup of maps when there are trajectories that touch or lie partially in the boundary.

Then for $t\geq 0$, 
\begin{equation}\label{expr phi-v}
\Phi^v_t(x)\ =\ x + \int_0^{t\wedge\tau_\partial^v(x)} v\bigl(\Phi^v_s(x)\bigr)\,ds,
\end{equation}
where $t\wedge\tau_\partial^v(x)$ denotes the minimum of $t$ and $\tau_\partial^v(x)$. The equivalent of Lemma \ref{lem: Pv mu - Pv' mu} can be obtained from \eqref{expr phi-v} by means of Gronwall's Lemma essentially, once one knows that the stopping time satisfies for fixed $x\in\overline{\Omega}$ a Lipschitz estimate of the form
\begin{equation}\label{essential estimate stopping time}
\bigl|t\wedge \tau^v_\partial(x) - t\wedge \tau^{v'}_\partial(x)\bigr|\leqslant C\|v-v'\|_\infty t.
\end{equation}
We are not aware of results in the literature that provide estimates like \eqref{essential estimate stopping time}. Neither did we succeed in establishing such an estimate ourselves. The rich possible dynamics of solutions of higher dimensional systems of non-linear differential equations may even make such global Lipschitz dependence of the velocity field impossible in general. Thus, a generalization of this part of the paper to two and higher dimensional state spaces seems not straightforward.
\end{remark}

\subsection{Properties of the solution for prescribed velocity}\label{sec: summ tech Q}
% General f (piecewise BL)
We consider mild solutions to \eqref{eq:main equation}, that are defined in the following sense:\\

\begin{definition}[See {\cite[Definition 2.4]{EvHiMu_JDE}}]\label{def: mild soln fixed v}
A {\em measure-valued mild solution} to the Cauchy-problem associated to \eqref{eq:main equation} on $[0,T]$ with initial value $\nu\in\CM([0,1])$ is a continuous map $\mu:[0,T]\to\CM([0,1])_\BL$ that is $\|\cdot\|_\TV$-bounded and that satisfies the variation of constants formula
\begin{equation}\label{eq:VOC}
\mu_t = P_t\,\nu \ +\ \int_0^t P_{t-s}F_f(\mu_s)\, ds\qquad \mbox{for all}\ t\in [0,T].
\end{equation}
Here, the perturbation map $\map{F_f}{\CM([0,1])}{\CM([0,1])}$ is given by $F_f(\mu):= f\cdot \mu$.
\end{definition}\\
\\
We showed in \cite{EvHiMu_JDE} that mild solutions in the sense of Definition~\ref{def: mild soln fixed v} exist, are unique and depend continuously on the initial data. We repeat those results in the following theorem.\\

\begin{theorem}\label{thm:exist uniq v prescribed}
Let $f:[0,1]\to\R$ be a piecewise bounded Lipschitz function such that $v(x)\neq 0$ at any point $x$ of discontinuity of $f$. Then for each $T\geqs 0$ and $\mu_0\in\CM([0,1])$ there exists a unique continuous and locally $\|\cdot\|_\TV$-bounded solution to \eqref{eq:VOC}.
Moreover, there exists $C_T>0$ such that for all initial values $\mu_0,\mu'_0\in\CM([0,1])$ the corresponding mild solutions $\mu$ and $\mu'$ satisfy
\begin{equation}\label{eq:cont dep initial cond}
\|\mu_t-\mu'_t\|_\BL^* \leqs C_T \|\mu_0-\mu'_0\|_\BL^*
\end{equation}
for all $t\in[0,T]$.
\end{theorem}\\

\begin{proof}
See {\cite[Propositions 3.1, 3.3 and 3.5]{EvHiMu_JDE}} for details.
\end{proof}\\
\\
In this paper, we restrict ourselves to those functions $f$ that are bounded Lipschitz on $[0,1]$; see \S \ref{sec: disc f} for further discussion on the need of this restriction. Let $v\in\BL([0,1])$ and $f\in\BL([0,1])$ be arbitrary. For all $t\geqs0$, we define $\map{Q_t}{\CM([0,1])}{\CM([0,1])}$ to be the operator that maps the initial condition to the solution in the sense of Definition \ref{def: mild soln fixed v}. Theorem \ref{thm:exist uniq v prescribed} guarantees that this operator is well-defined and continuous for $\|\cdot\|^*_\BL$. Moreover, $Q$ preserves positivity, due to \cite[Corollary 3.4]{EvHiMu_JDE}.\\
\\
In the rest of this section, we give an overview of the properties of the solution operator $Q$.\\

\begin{lemma}[Semigroup property]\label{lem: Q semigroup}
%For all $v\in\BL([0,1])$ and $f\in\BL([0,1])$, t
%The operator $Q_\cdot$
The set of operators $(Q_t)_{t\geqs0}$ satisfies the semigroup property. That is,
\begin{equation*}
Q_t\, Q_s\, \mu = Q_{t+s}\,\mu
\end{equation*}
for all $s,t\geqs0$ and for all $\mu\in\CM([0,1])$.
\end{lemma}\\

\begin{proof}
The proof follows the lines of argument of \cite[p.~283]{Sikic}. We consider
\begin{align}\label{eqn: Qts - QtQs written out}
\nonumber Q_{t+s}\mu - Q_t\,Q_s\,\mu =&\, P_{t+s}\mu + \Int{0}{t+s}{P_{t+s-\sigma}\,F_f(Q_\sigma\,\mu)}{d\sigma}\\
& - P_t\,Q_s\,\mu - \Int{0}{t}{P_{t-\sigma}\,F_f(Q_\sigma\,Q_s\,\mu)}{d\sigma},
\end{align}
and observe that
\begin{align}\label{eqn: PtQs}
\nonumber P_t\,Q_s\,\mu =&\, P_t\,P_s\,\mu + P_t\, \Int{0}{s}{P_{s-\sigma}\,F_f(Q_\sigma\,\mu)}{d\sigma}\\
=&\, P_{t+s}\,\mu + \Int{0}{s}{P_{t+s-\sigma}\,F_f(Q_\sigma\,\mu)}{d\sigma}.
\end{align}
Because $f\in\BL([0,1])$, the map $\sigma\mapsto P_{s-\sigma}\,F_f(Q_\sigma\,\mu)$ is continuous and hence it is measurable. Therefore, the second equality in \eqref{eqn: PtQs} holds due to \cite[Equation~(C.3)]{EvHiMu_JDE}. %\eqref{eqn: semigroup inside integral}. %{\color{red} measurability integrand, ref. Ch 4   Remark \ref{rem: pert Bochner meas}}%the integrand is even continuous for BL f.  %the integrals are well-defined due to the theorems in Ch 4
A combination of \eqref{eqn: Qts - QtQs written out} and \eqref{eqn: PtQs} yields that
\begin{align}
\nonumber Q_{t+s}\mu - Q_t\,Q_s\,\mu =& \Int{s}{t+s}{P_{t+s-\sigma}\,F_f(Q_\sigma\,\mu)}{d\sigma} - \Int{0}{t}{P_{t-\sigma}\,F_f(Q_\sigma\,Q_s\,\mu)}{d\sigma}\\
=& \Int{0}{t}{P_{t-\sigma}\,\left(F_f(Q_{\sigma+s}\,\mu)-  F_f(Q_\sigma\,Q_s\,\mu)\right)}{d\sigma}.\label{eqn: Qts - QtQs = integral}
\end{align}
To obtain the last step in \eqref{eqn: Qts - QtQs = integral}, we use the coordinate transformation $\tau:=\sigma-s$ in the first integral and subsequently renamed the new variable $\tau$ as $\sigma$. We estimate the total variation norm of \eqref{eqn: Qts - QtQs = integral} in the following way:
\begin{align*}
\nonumber \|Q_{t+s}\mu - Q_t\,Q_s\,\mu\|_\TV \leqs& \Int{0}{t}{\|P_{t-\sigma}\,\left(F_f(Q_{\sigma+s}\,\mu)-  F_f(Q_\sigma\,Q_s\,\mu)\right)\|_\TV}{d\sigma}\\
\nonumber \leqs& \Int{0}{t}{\|F_f(Q_{\sigma+s}\,\mu)-  F_f(Q_\sigma\,Q_s\,\mu)\|_\TV}{d\sigma}\\
\leqs& \, \|f\|_\infty\,\Int{0}{t}{\|Q_{\sigma+s}\,\mu- Q_\sigma\,Q_s\,\mu\|_\TV}{d\sigma}.
\end{align*}
Here, we used \cite[Proposition C.2(iii)]{EvHiMu_JDE} %Part \ref{prop: part TV norm inside integral} of Proposition \ref{prop:TV-estimates} 
(noting that the integrands are continuous with respect to $\sigma$) in the first line, \eqref{eqn:TV norm Pt} in the second line and the fact that $f\in\BL([0,1])\subset C_b([0,1])$ in the last line. Gronwall's Lemma now implies that $\|Q_{t+s}\mu - Q_t\,Q_s\,\mu\|_\TV=0$ for all $s,t\geqs0$.
\end{proof}\\

\begin{lemma}\label{lem: Qtmu - Qsmu}
For all $\mu\in\CM([0,1])$ and $s,t\geqs0$, we have that
\[\|Q_t\mu - Q_s\mu\|^*_\BL \leqs \|\mu\|_\TV\cdot \bigl(\|f\|_\infty + \|v\|_\infty\bigr)\cdot e^{\|f\|_\infty\max(t,s)}\cdot |t-s|.\]
\end{lemma}\\

\begin{proof}
The statement of this lemma is part of the result of \cite[Proposition 3.3]{EvHiMu_JDE}.%Proposition \ref{prop:existence}.
\end{proof}\\

\begin{lemma}\label{lem: TV BL bounds Qmu}
%For all $v\in\BL([0,1])$, $f\in\BL([0,1])$,
For all $\mu\in\CM([0,1])$ and $t\geqs0$, we have that
\begin{enumerate}
\item[(i)] $\| Q_t\mu\|_\TV \leqs \|\mu\|_\TV \, \exp(\|f\|_\infty\,t)$, and  \label{lem: part TV bound Qmu}
\item[(ii)] $\| Q_t\mu\|_\BL^* \leqs \|\mu\|^*_\BL\, \exp(|v|_\Lip\,t+\|f\|_\BL\,t\,e^{|v|_\Lip\,t})$.\label{lem: part BL bound Qmu}
\end{enumerate}
\end{lemma}
\vspace{\baselineskip}
\begin{proof}
(i): This estimate is given in \cite[Proposition 3.3]{EvHiMu_JDE}.\\%Proposition \ref{prop:existence}.\\
(ii): By applying \cite[(C.1)]{EvHiMu_JDE} %\eqref{eqn: Bochner norm inside} 
and Lemma \ref{lemma:lift is Lipschitz in time}(ii) %\ref{lem: part BL bound Pmu} 
we obtain from \eqref{eq:VOC} the estimate
\[\| Q_t\,\mu\|_\BL^* \leqs \exp(|v|_\Lip\,t)\,\|\mu\|^*_\BL + \Int{0}{t}{\exp(|v|_\Lip\,(t-s))\|f\|_\BL\,\|Q_s\mu\|^*_\BL}{ds}.  \]
Gronwall's Lemma now yields the statement of Part (ii) %\ref{lem: part BL bound Qmu} 
of the lemma.
\end{proof}\\

\begin{corollary}\label{cor: Qmu - Qnu}
For all $\mu,\nu\in\CM([0,1])$ and $t\geqs0$, we have that
\[\|Q_t\,\mu - Q_t\,\nu\|^*_\BL \leqs \|\mu-\nu\|^*_\BL\, \exp(|v|_\Lip\,t+\|f\|_\BL\,t\,e^{|v|_\Lip\,t}).\]
\end{corollary}\\

\begin{proof}
Apply Part (ii) %\ref{lem: part BL bound Qmu} 
of Lemma \ref{lem: TV BL bounds Qmu} to the measure $\mu-\nu\in\CM([0,1])$.
\end{proof}\\
\\
%Similar to the semigroups $P^v$ and $P^{v'}$, w
We write $Q^v$ and $Q^{v'}$ to distinguish between the semigroups $Q$ on $\CM([0,1])$ associated to $v\in\BL([0,1])$  and $v'\in\BL([0,1])$, respectively.\\

\begin{lemma}\label{lem: Qvmu - Qv'mu}
For all $v,v'\in\BL([0,1])$, $\mu\in\CM([0,1])$ and $t\geqs0$, the following estimate holds:
\[\|Q^v_t\mu - Q^{v'}_t\mu\|^*_\BL \leqs \|v-v'\|_\infty\,\|\mu\|_\TV\,\exp(L\,t+\|f\|_\BL\,t\,e^{L\,t})\cdot[t+t^2\|f\|_\infty \,e^{\|f\|_\infty\,t}],\]
where $L:=\min(|v|_\Lip\,,|v'|_\Lip)$.
\end{lemma}\\

\begin{proof}
We have
\begin{equation}\label{eqn: Qv - Qv'}
\|Q^v_t\mu - Q^{v'}_t\mu\|^*_\BL \leqs \|P^v_t\mu - P^{v'}_t\mu\|^*_\BL + \Int{0}{t}{\| P^v_{t-s}F_f(Q^v_s\mu) - P^{v'}_{t-s}F_f(Q^{v'}_s\mu) \|_\BL^*}{ds}.
\end{equation}
Lemma \ref{lem: Pv mu - Pv' mu} provides an appropriate estimate of the first term on the right-hand side. For the integrand in the second term, we have
\begin{align}
\nonumber \| P^v_{t-s}F_f(Q^v_s\mu) - P^{v'}_{t-s}F_f(Q^{v'}_s\mu) \|_\BL^* \leqs&\, \| P^v_{t-s}F_f(Q^v_s\mu) - P^{v'}_{t-s}F_f(Q^{v}_s\mu) \|_\BL^*\\
\nonumber & + \| P^{v'}_{t-s}F_f(Q^v_s\mu) - P^{v'}_{t-s}F_f(Q^{v'}_s\mu) \|_\BL^*\\
\nonumber \leqs&\, \|v-v'\|_\infty\,(t-s)\,\|F_f(Q^v_s\mu)\|_\TV\,e^{L(t-s)}\\
& + e^{|v'|_\Lip\,(t-s)}\,\| F_f(Q^v_s\mu) - F_f(Q^{v'}_s\mu) \|^*_\BL,\label{eqn: PvFQv - Pv'FQv' first step}
\end{align}
due to Lemma \ref{lem: Pv mu - Pv' mu} and Lemma \ref{lemma:lift is Lipschitz in time}(ii). %\ref{lem: part BL bound Pmu}. 
We proceed by estimating the right-hand side of \eqref{eqn: PvFQv - Pv'FQv' first step} and obtain
\begin{align}
\nonumber \| P^v_{t-s}F_f(Q^v_s\mu) - P^{v'}_{t-s}F_f(Q^{v'}_s\mu) \|_\BL^* \leqs&\,  \|v-v'\|_\infty\,(t-s)\,\|f\|_\infty\,\|\mu\|_\TV\,e^{\|f\|_\infty\,s}\,e^{L(t-s)}\\
& + e^{|v'|_\Lip\,(t-s)}\,\|f\|_\BL\,\| Q^v_s\mu - Q^{v'}_s\mu \|^*_\BL,\label{eqn: PvFQv - Pv'FQv'}
\end{align}
where we use Part (i) %\ref{lem: part TV bound Qmu} 
of Lemma \ref{lem: TV BL bounds Qmu} in the first term on the right-hand side. Since the estimate in \eqref{eqn: PvFQv - Pv'FQv'} is symmetric in $v$ and $v'$, we can replace $|v'|_\Lip$ by $L$.\\
Substitution of the result of Lemma \ref{lem: Pv mu - Pv' mu} and \eqref{eqn: PvFQv - Pv'FQv'} in \eqref{eqn: Qv - Qv'} yields
\begin{align*}
\nonumber \|Q^v_t\mu - Q^{v'}_t\mu\|^*_\BL \leqs&\,  \|v-v'\|_\infty\,t\,\|\mu\|_\TV\,e^{Lt}\,(1+t\,\|f\|_\infty\,\,e^{\|f\|_\infty\,t})\\
& + e^{Lt}\,\|f\|_\BL\,\Int{0}{t}{\| Q^v_s\mu - Q^{v'}_s\mu \|^*_\BL}{ds}.
\end{align*}
The statement of the lemma follows from Gronwall's Lemma.
\end{proof}

\section{Measure-dependent velocity fields: main results}\label{sec: gen meas-dep v}
This section contains the main results of the present work. We generalize the assumptions on $v$ from \cite{EvHiMu_JDE} in the following way to measure-dependent velocity fields:\\

\begin{assumption}[Assumptions on the measure-dependent velocity field]\label{ass: v properties}
Assume that $\map{v}{\CM([0,1])\times[0,1]}{\R}$ is a mapping such that:
\begin{enumerate}
  \item[(i)] $v[\mu]\in\BL([0,1])$, for each $\mu\in\CM([0,1])$.\label{ass: part v meas-dep in BL}
\end{enumerate}
Furthermore, assume that for any $R>0$ there are constants $K_R$, $L_R$, $M_R$ such that for all $\mu,\nu\in\CM([0,1])$ satisfying $\|\mu\|_\TV\leqslant R$ and $\|\nu\|_\TV\leqslant R$, the following estimates hold:
\begin{enumerate}
  \item[(ii)] $\|v[\mu]\|_\infty\leqs K_R$,\label{ass: part unif bnd}
  \item[(iii)] $|\,v[\mu]\,|_\Lip\leqs L_R$, \,\,\, and\label{ass: part Lipsch in x}
  \item[(iv)] $\|v[\mu]-v[\nu]\|_\infty \leqslant M_R\,\|\mu-\nu\|^*_\BL$.\label{ass: part Lipsch in mu incl R}
\end{enumerate}
\end{assumption}
\vspace{\baselineskip}
\begin{example}\label{ex: meas-dep v}
An example of a function $v$ satisfying Assumption \ref{ass: v properties} is:
\begin{equation}\label{eqn: example meas-dep v}
v[\mu](x) := \Int{[0,1]}{}{\K(x-y)}{d\mu(y)} = (\K*\mu)(x),
\end{equation}
for each $\mu\in\CM([0,1])$ and $x\in[0,1]$, with $\K\in\BL([-1,1])$. This is a relevant choice, because it models interactions among individuals.
\end{example}\\

\begin{remark}\label{rem:estimates on ball i.o. unif Hoogwater}
Parts (ii) and (iii) of Assumption \ref{ass: v properties} are an improvement compared to \cite{Hoogwater-thesis}. There, the infinity norm and Lipschitz constant are assumed to hold uniformly for all $\mu\in\CM([0,1])$; cf.~Assumption (F1) on \cite[p.~40]{Hoogwater-thesis}. We note that the convolution in Example \ref{ex: meas-dep v} satisfies Assumption \ref{ass: v properties}, but does not satisfy Assumption (F1) in \cite{Hoogwater-thesis}. They require a uniform Lipschitz constant because their Lemma 4.3 is an estimate in the $\|\cdot\|^*_\BL$-norm for which Part (ii) %\ref{lem: part BL bound Pmu} 
of our Lemma \ref{lemma:lift is Lipschitz in time} is used. Our counterpart of Lemma 4.3 in \cite{Hoogwater-thesis} is Lemma \ref{lem: set timeslices bdd}. We give an estimate in terms of the $\|\cdot\|_\TV$-norm using \eqref{eqn:TV norm Pt} which does not involve the Lipschitz constant.
\end{remark}\\
\\
%Let $v$ be a mapping from $\CM([0,1])$ to the space of functions $[0,1]\to\R$. Let $\mu\in\CM([0,1])$ be fixed but arbitrary. We assume (in agreement with the assumptions on $v$ in \S \ref{sec: summ JDE paper}) that $v[\mu]$ is bounded Lipschitz on $(0,1)$. In agreement with the assumptions on $v$ in \S \ref{sec: summ JDE paper}, we assume that $v[\mu](x)=0$ at $x=0$ or $x=1$, whenever $\lim_{x\downarrow0}v[\mu](x)<0$ or $\lim_{x\uparrow1}v[\mu](x)>0$, respectively.\\    !!!!!!!: Section heading \ref{sec: summ JDE paper} was removed!!
%\\
Our aim is to prove well-posedness (in some sense yet to be defined) of \eqref{eqn:main equation nonlin}. That is,
\begin{equation*}
\frac{\partial}{\partial t}\mu_t + \frac{\partial}{\partial x} (v[\mu_t]\,\mu_t) = f\cdot \mu_t
\end{equation*}
on $[0,1]$. %where we focus on the physically relevant \textit{positive} measure-valued solutions (cf.~Remark \ref{rem: pos meas interpretable}).
As said in \S \ref{sec: summ tech Q}, we restrict ourselves to $f$ that is bounded Lipschitz on $[0,1]$.\\
\\
%where we restrict ourselves to $f$ that is bounded Lipschitz on $[0,1]$; cf.~\S \ref{sec: discussion bc meas-dep v} for further discussion on this restriction. We focus on the physically relevant \textit{positive} measure-valued solutions (cf.~Remark \ref{rem: pos meas interpretable}), hence $v$ was defined only for $\mu\in\CM^+([0,1])$.\\
%\\
%\textbf{\color{red}Note that for the $T/2^k$ spacing, we start counting at $j=1$!}
%Let $T>0$ be given. Let $(d_k)_{k\in\N}\subset\Rp$ be a strictly decreasing sequence such that $\lim_{k\to\infty}d_k=0$. For each $k\in\N$, let $\underline{t}^k:=(t_1^k,\ldots,t_{N_k}^k)^\top\in[0,T]^{N_k}$ be a partition of $[0,T]$. Here, $(N_k)_{k\in\N}\subset\N$ is a strictly increasing sequence. For each $k\in\N$ we require that $t_1^k=0$, $t_{N_k}^k=T$, $t_i^k<t_{i+1}^k$ for each $i\in\{1,\ldots,N_k-1\}$, and $\max_{i\in\{1,\ldots,N_k-1\}}t_{i+1}^k-t_i^k\leqslant d_k$. Moreover, we require that $\underline{t}^{k+1}$ is a \textit{refinement} of $\underline{t}^k$: i.e.~for all $i\in\{1,\ldots,N_k\}$ there is a $j\in\{1,\ldots,N_{k+1}\}$ such that $t_i^k=t_j^{k+1}$.\\
%\\
We now introduce the aforementioned forward-Euler-like approach to construct approximate solutions. Let $T>0$ be given. Let $N\geqslant 1$ be fixed and define a set $\alpha\subset [0,T]$ as follows: 
\begin{equation}\label{eqn: partition}
\alpha:= \big\{ t_j\in[0,T] \;\, : \;\, 0\leqslant j \leqslant N,\; t_0=0,\; t_{N}=T,\; t_j<t_{j+1}  \big\}.
\end{equation}
A set $\alpha$ of this form is called a \textit{partition} of the interval $[0,T]$ and $N$ denotes the number of \textit{subintervals} in $\alpha$.\\
\\
Let $\mu_0\in\CM([0,1])$ be fixed. For a given partition $\alpha:=\{t_0,\ldots,t_N\}\subset[0,T]$, define a measure-valued trajectory $\mu \in C([0,T];\CM([0,1]))$ by %Separate definition k=0 (deviating from k>0) not necessary
\begin{equation}\label{eqn: Euler scheme}
\left\{
  \begin{array}{ll}
    \mu_t := Q^{v_j}_{t-t_j}\,\mu_{t_j},&\text{if }t\in(t_j,t_{j+1}];\vspace{0.2cm} \\ %P^{v_j}_{t-t_j}\mu_{t_j}+\ \displaystyle\int_{t_j}^t P^{v_j}_{t-s}F_f(\mu_s)\, ds,&\text{if }t\in(t_j,t_{j+1}];\\
    v_j:=v[\mu_{t_j}];&\vspace{0.2cm}\\
    \mu_{t=0}=\mu_0,&
  \end{array}
\right.
\end{equation}
for all $j\in\{0,\ldots,N-1\}$. Here, $(Q^v_t)_{t\geqslant0}$ denotes the semigroup introduced in \S \ref{sec: summ tech Q} associated to an arbitrary $v\in\BL([0,1])$. Note that by Assumption \ref{ass: v properties}, Part (i), $v_j=v[\mu_{t_j}]\in\BL([0,1])$ for each $j$.\\     %of the stopped flow associated to any $v\in\BL([0,1])$ as given by \eqref{Def Pt push forward}.\\
\\
We call this a forward-Euler-like approach, because it is the analogon of the forward Euler method for ODEs (cf.~e.g.~\cite[Chapter 2]{Butcher}). Consider the %toy
ODE $dx/dt = v(x)$ on $\R$ for some (Lipschitz continuous) $\map{v}{\R}{\R}$. The forward Euler method approximates the solution on some interval $(t_j, t_{j+1}]$ by evolving the approximate solution at time $t_j$, named $x_j$, due to a \textit{constant} velocity $v(x_j)$. That is, $x(t)\approx x_j + (t-t_j)\cdot v(x_j)$ for all $t\in (t_j, t_{j+1}]$.\\
In \eqref{eqn: Euler scheme}, we introduce the approximation $\mu_t$, where $\mu_t$ results from $\mu_{t_j}$ by the evolution due to the constant velocity field $v[\mu_{t_j}]$. The word \textit{constant} here does not refer to $v$ being the same for all $x\in[0,1]$, but to the fact that $v$ corresponding to the same $\mu_{t_j}$ is used throughout $(t_j, t_{j+1}]$.\\% Equation \eqref{eqn: Euler scheme} shows that actually obtaining $\mu_t^k$ is not as straightforward as for the ODE mentioned above.\\
\\
The conditions in Parts (ii)--(iv) %\ref{ass: part unif bnd}--\ref{ass: part Lipsch in mu incl R} 
of Assumption \ref{ass: v properties} are only required to hold for measures in a $\TV$-norm bounded set, in view of the following lemma:\\

\begin{lemma}\label{lem: set timeslices bdd}
Let $\mu_0\in\CM([0,1])$ be given and let $\map{v}{\CM([0,1])\times[0,1]}{\R}$ satisfy Assumption \ref{ass: v properties}(i). %\ref{ass: part v meas-dep in BL}. 
For a given partition $\alpha:=\{t_0,\ldots,t_N\}\subset[0,T]$, let $\mu \in C([0,T];\CM([0,1]))$ be defined by \eqref{eqn: Euler scheme}. Then the set of all \textit{timeslices} of $\mu$, that is
\begin{equation*}
\mathcal{A}:= \{ \mu_t: t\in[0,T] \},
\end{equation*}
is bounded in both $\|\cdot\|_\TV$ and $\|\cdot\|^*_\BL$. The bounds are independent of the choice of $\alpha$.
\end{lemma}\\

\begin{proof}
Fix $j\in\{0,\ldots,N-1\}$ and let $t\in(t_j,t_{j+1}]$. By Part (i) %\ref{lem: part TV bound Qmu} 
of Lemma \ref{lem: TV BL bounds Qmu}, we have that
\begin{align*}
\nonumber \|\mu_t\|_\TV = \|Q^{v_j}_{t-t_j}\,\mu_{t_j}\|_\TV \leqs& \,\|\mu_{t_j}\|_\TV \,\exp(\|f\|_\infty\,(t-t_j))\\
\leqs& \, \|\mu_{t_j}\|_\TV \, \exp(\|f\|_\infty\,(t_{j+1}-t_j))
\end{align*}
%Then we have that
%\begin{align}
%\nonumber \|\mu^k_t\|_\TV &\leqslant \|P^{v^k_j}_{t-t^k_j}\mu^k_{t^k_j}\|_\TV +\ \displaystyle\int_{t^k_j}^t \|P^{v^k_j}_{t-s}F_f(\mu_s^k)\|_\TV\, ds\\
%&\leqslant \|\mu^k_{t^k_j}\|_\TV + \|f\|_\infty\,\displaystyle\int_{t^k_j}^t \|\mu_s^k\|_\TV\, ds,
%\end{align}
%due to Part \ref{prop: part TV norm inside integral} of Proposition \ref{prop:TV-estimates}, \eqref{eqn:TV norm Pt} and the fact that $f\in\BL([0,1])\subset C_b([0,1])$.
%%where we have used \eqref{eqn:TV norm Pt} %Part \ref{lem: part BL bound Pmu} of Lemma \ref{lemma:lift is Lipschitz in time}
%%and \eqref{eqn: Bochner norm inside}.
%We use Gronwall's Lemma and estimate the length of the domain of integration $t-t^k_j$ against $T/2^k$. We obtain
%\begin{equation}
%\|\mu^k_t\|_\TV \leqslant \|\mu^k_{t^k_j}\|_\TV\,\exp(\|f\|_\infty\,T/2^k),
%\end{equation}
for all $t\in(t_j,t_{j+1}]$. Iteration of the right-hand side with respect to $j$ yields
\begin{equation*}
\|\mu_t\|_\TV \;\leqslant\; \|\mu_0\|_\TV\,\prod_{i=0}^j \exp(\|f\|_\infty\,(t_{i+1}-t_i)) \;=\; \|\mu_0\|_\TV\,\exp(\|f\|_\infty\,(t_{j+1}-t_0)).
\end{equation*}
Hence, for all $t\in[0,T]$
\begin{equation*}
\|\mu_t\|_\TV \;\leqslant\; \|\mu_0\|_\TV \,\exp(\|f\|_\infty\,(t_{N}-t_0)) \;=\; \|\mu_0\|_\TV \,\exp(\|f\|_\infty\,T).
\end{equation*}
This bound is in particular independent of $t$, $N$ and the distribution of points within $\alpha$. The bound in $\|\cdot\|^*_\BL$ follows from the inequality $\|\nu\|^*_\BL\leqs \|\nu\|_\TV$ that holds for all $\nu\in\CM([0,1])$.
%where the constant $A$ is independent of $k$ and $t$. Since positivity is preserved along solutions of \eqref{eqn: Euler scheme}, $\mu^k_t$ is a positive measure for any $k$ and $t$, and hence $\|\mu^k_t\|_\TV =\|\mu^k_t\|^*_\BL\leqslant A\,\|\mu_0\|^*_\BL$ is bounded.
\end{proof}\\
\\
In this paper we construct sequences of Euler approximations, each following from a sequence of partitions $(\alpha_k)_{k\in\N}$ that satisfies the following assumption:
\begin{assumption}[Assumptions on the sequence of partitions]\label{ass: partition}
Let $(\alpha_k)_{k\in\N}$ be a sequence of partitions of $[0,T]$ and let $(N_k)_{k\in\N}\subset\N$ be the corresponding sequence such that each $\alpha_k$ is of the form
\begin{equation}\label{eqn: partition alpha_k}
\alpha_k:= \big\{ t^k_j\in[0,T] \;\, : \;\, 0\leqslant j \leqslant N_k,\; t^k_0=0,\; t^k_{N_k}=T,\; t^k_j<t^k_{j+1}  \big\}.
\end{equation}
Define 
\begin{equation}
M^{(k)}:= \max_{j\in\{0,\ldots,N_k-1\}}t^k_{j+1}-t^k_j \label{eqn: def max interval}
\end{equation}
for all $k\in\N$. Assume that the sequence $(M^{(k)})_{k\in\N}$ is nonincreasing and $M^{(k)}\to 0$ as $k\to \infty$.
%\begin{enumerate}
%  \item[(i)] the sequence $(N_k)_{k\in\N}$ is nondecreasing; {\color{red}\bf NEEDED?!}
%  \item[(ii)] the sequence $(M^{(k)})_{k\in\N}$ is nonincreasing and $M^{(k)}\to 0$ as $k\to \infty$.
%\end{enumerate}
\end{assumption}\\

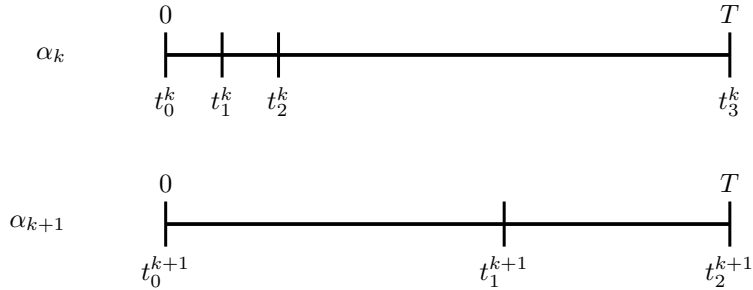
\begin{figure}[h]
\begin{center}
        \begin{tikzpicture}[scale=1.5, >= latex]
        \begin{scope}[shift={(0,0)}]
        \draw (-0.8,0)node[anchor=east]{$\alpha_k$};
        \draw[-, line width = 0.5mm] (0,0)--(5,0);
        \draw[-, line width = 0.4mm] (0,0.2)--(0,-0.2);
        \draw (0,0.2)node[anchor=south]{$0$};
        \draw (0,-0.2)node[anchor=north]{$t^k_0$};
        \draw[-, line width = 0.4mm] (5,0.2)--(5,-0.2);
        \draw (5,0.2)node[anchor=south]{$T$};
        \draw (5,-0.2)node[anchor=north]{$t^k_3$};
        
        \draw[-, line width = 0.4mm] (0.5,0.2)--(0.5,-0.2);
        \draw (0.5,-0.2)node[anchor=north]{$t^k_1$};
        \draw[-, line width = 0.4mm] (1,0.2)--(1,-0.2);
        \draw (1,-0.2)node[anchor=north]{$t^k_2$};
        \end{scope}

		\begin{scope}[shift={(0,-1.5)}]
		\draw (-0.8,0)node[anchor=east]{$\alpha_{k+1}$};
        \draw[-, line width = 0.5mm] (0,0)--(5,0);
        \draw[-, line width = 0.4mm] (0,0.2)--(0,-0.2);
        \draw (0,0.2)node[anchor=south]{$0$};
        \draw (0,-0.2)node[anchor=north]{$t^{k+1}_0$};
        \draw[-, line width = 0.4mm] (5,0.2)--(5,-0.2);
        \draw (5,0.2)node[anchor=south]{$T$};
        \draw (5,-0.2)node[anchor=north]{$t^{k+1}_2$};
        
        \draw[-, line width = 0.4mm] (3,0.2)--(3,-0.2);
        \draw (3,-0.2)node[anchor=north]{$t^{k+1}_1$};
		\end{scope}
        \end{tikzpicture}
\caption{Two possible subsequent partitions in a sequence $(\alpha_k)_{k\in\N}$ satisfying Assumption \ref{ass: partition}.}\label{fig: partitions}
\end{center}
\end{figure}

\begin{example}\label{ex: partitions}
The following sequences of partitions satisfy Assumption \ref{ass: partition}:
\begin{itemize}
\item  For all $k\in\N$, take $N_k:=2^k$, and let $t^k_j:= jT/2^k$ for all $j\in\{0,\ldots,N_k\}$. This implies that $M^{(k)}=T/2^k$ for all $k\in\N$. This specific sequence of partitions was used in \cite[Chapter 5]{Evers_PhD}.
\item Fix $q\in\Np$. For all $k\in\N$, take $N_k:=q^k$, and let $t^k_j:= jT/q^k$ for all $j\in\{0,\ldots,N_k\}$. This implies that $M^{(k)}=T/q^k$ for all $k\in\N$. In the discussion section of \cite[Chapter 5]{Evers_PhD}, the results of the current paper were conjectured to hold for this case. 
\item For all $k\in\N$, take $N_k:=k+1$, and let $t^k_j:= jT/(k+1)$ for all $j\in\{0,\ldots,N_k\}$. This implies that $M^{(k)}=T/(k+1)$ for all $k\in\N$. This is an elementary time discretization (with uniform mesh size) used frequently when proving the convergence of numerical methods. 
\item Let $\alpha_0$ be a possibly non-uniform partition of $[0,T]$. Construct the sequence $(\alpha_k)_{k\in\N}$ in such a way that any $\alpha_{k+1}$ is a refinement of $\alpha_k$.  That is, $\alpha_{k+1}\subset\alpha_k$ for all $k\in\N$. Elements may be added in a non-uniform fashion to obtain $\alpha_{k+1}$ from $\alpha_k$, as long as $M^{(k)}\to 0$ as $k\to\infty$. In this case $(N_k)_{k\in\N}$ is automatically nondecreasing.
\end{itemize}
Also, some less straightforward sequences of non-uniform partitions are admissible, in which subsequent partitions are not refinements. See for example Figure \ref{fig: partitions}, in which two subsequent elements from the sequence $(\alpha_k)_{k\in\N}$ are given. These elements could indeed occur, since $M^{(k+1)}<M^{(k)}$. This example is rather counter-intuitive, as there is a local growth of the mesh size at the left-hand side of the interval $[0,T]$ when we go from $\alpha_k$ to $\alpha_{k+1}$. Note that even $N^{(k+1)}<N^{(k)}$. However, admissibility of a sequence of partitions is only determined by the local ordering of the \emph{maximum} mesh spacing (i.e. the condition $M^{(k+1)}\leqslant M^{(k)}$) and its long-time behaviour: $M^{(k)}\to 0$ as $k\to\infty$.
\end{example}\\

\begin{remark}
Assumption \ref{ass: partition} implies that $N_k\to\infty$ as $k\to\infty$.\\
If $(M^{(k)})_{k\in\N}$ is \emph{not} nonincreasing, but still $M^{(k)}\to 0$, then it is possible to extract a subsequence $(\alpha_{k_\ell})_{\ell\in\N}$ such that $(M^{(k_\ell)})_{\ell\in\N}$ \emph{is} nonincreasing.
\end{remark}\\
\\
We define a mild solution in this context as follows:\\ %-- cf.~\cite{Gwiazda1}
\begin{definition}[Mild solution of \eqref{eqn:main equation nonlin}]\label{def: mild soln nonlin}
Let the space of continuous maps from $[0,T]$ to $\CM([0,1])_\BL$ be endowed with the metric defined for all $\mu,\nu\in C([0,T]; \CM([0,1]))$ by
\begin{equation}\label{eqn: metric sup BL norm}
\sup_{t\in[0,T]}\|\mu_t-\nu_t\|^*_\BL.
\end{equation}
Let $(\alpha_k)_{k\in\N}$ be a sequence of partitions satisfying Assumption \ref{ass: partition}. For each $k\in\N$, let $\mu^k\in C([0,T];\CM([0,1]))$ be defined by \eqref{eqn: Euler scheme} with partition $\alpha_k$. Then, for any such sequence of partitions $(\alpha_k)_{k\in\N}$, any limit of a subsequence of $(\mu^k)_{k\in\N}$ is called a (measure-valued) mild solution of \eqref{eqn:main equation nonlin}.
\end{definition}\\
\\
The name \textit{mild solutions} is appropriate, because they are constructed from piecewise mild solutions in the sense of Definition \ref{def: mild soln fixed v}.\\

\begin{remark}
Consider the solution of \eqref{eqn: Euler scheme} for any partition $\alpha\subset[0,T]$. Mass that has accumulated on the boundary can move back into the interior of the domain whenever the velocity changes direction from one time interval to the next. This is due to the definition of the maximal interval of existence $I_{x_0}$ and the hitting time $\tau_\partial(x_0)$ in \S \ref{sec: prop stopped flow}. %\ref{sec:inidividualistic flow}. %Moving back can thus also happen in the limit solution
\end{remark}\\
\\
%\textbf{\color{red} This definition seems to be inspired by some kind of compactness of the function space. Can we say more about this?\\
%This def is also used in \cite{Gwiazda1} (?)}\\
%\\
In the rest of this paper we focus on \textit{positive} measure-valued solutions, because these are the only physically relevant solutions in many applications. %(cf.~Remark \ref{rem: pos meas interpretable}). 
The main result of this paper is the following theorem.\\

\begin{theorem}\label{thm: exist uniq nonlin}
Let $\mu_0\in\CM^+([0,1])$ be given and let $\map{v}{\CM([0,1])\times[0,1]}{\R}$ satisfy Assumption \ref{ass: v properties}. Endow the space $C([0,T];\CM([0,1]))$ with the metric defined by \eqref{eqn: metric sup BL norm}. Then, there is a unique element of $C([0,T]; \CM^+([0,1]))$ with initial condition $\mu_0$, that is a mild solution in the sense of Definition \ref{def: mild soln nonlin}. That is, for each sequence of partitions $(\alpha_k)_{k\in\N}$ satisfying Assumption \ref{ass: partition}, the corresponding sequence $(\mu^k)_{k\in\N}$ defined by \eqref{eqn: Euler scheme} is a sequence in $C([0,T];\CM^+([0,1]))$ and has a unique limit as $k\to\infty$.\\
Moreover, this limit is independent of the choice of $(\alpha_k)_{k\in\N}$.
\end{theorem}\\

%\begin{theorem}\label{thm: exist uniq nonlin}
%Let $\mu_0\in\CM^+([0,1])$ be given and let $\map{v}{\CM([0,1])\times[0,1]}{\R}$ satisfy Assumption \ref{ass: v properties}. Endow the space $C([0,T];\CM([0,1]))$ with the metric defined by \eqref{eqn: metric sup BL norm}. Then, for each sequence of partitions $(\alpha_k)_{k\in\N}$ satisfying Assumption \ref{ass: partition}, the corresponding sequence $(\mu^k)_{k\in\N}$ defined by \eqref{eqn: Euler scheme} is a sequence in $C([0,T];\CM^+([0,1]))$ and has a unique limit as $k\to\infty$. 
%\end{theorem}\\
%
%\begin{theorem}\label{thm: uniq indep of partition}
%For each $\mu_0\in\CM^+([0,1])$ and $\map{v}{\CM([0,1])\times[0,1]}{\R}$ satisfying Assumption \ref{ass: v properties}, there is a unique element of $C([0,T]; \CM^+([0,1]))$ with initial condition $\mu_0$, that is a mild solution in the sense of Definition \ref{def: mild soln nonlin}. That is, the unique limit in Theorem \ref{thm: exist uniq nonlin} is independent of the choice of $(\alpha_k)_{k\in\N}$, provided that Assumption \ref{ass: partition} is satisfied.
%\end{theorem}\\

\begin{corollary}[Global existence and uniqueness]\label{cor: global exist uniq}
For each $\mu_0\in\CM^+([0,1])$ and $\map{v}{\CM([0,1])\times[0,1]}{\R}$ satisfying Assumption \ref{ass: v properties}, a unique mild solution exists for all time $t\geqslant0$.
\end{corollary}\\

\begin{theorem}[Continuous dependence on initial data]\label{thm: cont dep}
For all $T>0$ and $\tilde{R}>0$ there is a constant $C_{\tilde{R},T}$ such that for all $\mu_0,\nu_0\in\CM^+([0,1])$ satisfying $\|\mu_0\|_\TV \leqs \tilde{R}$ and $\|\nu_0\|_\TV\leqs \tilde{R}$, the corresponding mild solutions $\mu,\nu\in C([0,T]; \CM^+([0,1]))$ satisfy
\begin{equation*}
\sup_{\tau\in[0,T]}\|\mu_\tau-\nu_\tau\|^*_\BL \leqslant C_{\tilde{R},T}\,\|\mu_0-\nu_0\|^*_\BL.
\end{equation*}
\end{theorem}\\
\\
The proofs of these theorems and this corollary are given in the next section, \S \ref{sec: proofs ch bc nonlin}. The key idea of the proof of Theorem \ref{thm: exist uniq nonlin} is to show that the sequence $(\mu^k)_{k\in\N}$ is a Cauchy sequence in a complete metric space, hence converges. We use estimates between approximations $\mu^k$ and $\mu^m$, $m\geqslant k$. Similar estimates are employed to obtain the result of Theorem \ref{thm: cont dep}. %The result of Theorem \ref{thm: uniq indep of partition} follows from an argument where we take two sequences of partitions $(\alpha_k)$ and $(\beta_k)$ and construct a third one $(\gamma_k)$ such that it contains subsequences of $(\alpha_k)$ and $(\beta_k)$. We show that the corresponding limit solutions coincide with the limit associated to $(\alpha_k)$ as well as the limit associated to $(\beta_k)$. 
To prove Corollary \ref{cor: global exist uniq}, we show that a solution at time $t\geqslant 0$ is provided by Theorem \ref{thm: exist uniq nonlin}, if $T>0$ is chosen such that $t\in[0,T]$. Moreover, this solution at time $t$ is independent of the exact choice of $T$.

%{\color{red}Structure/philosophy of the proof}

\section{Proofs of Theorems \ref{thm: exist uniq nonlin} and \ref{thm: cont dep}, and of Corollary \ref{cor: global exist uniq}}\label{sec: proofs ch bc nonlin}
In this section we prove the main results of this paper: Theorem \ref{thm: exist uniq nonlin}, Corollary \ref{cor: global exist uniq} and Theorem \ref{thm: cont dep}. The essential part of the proof of Theorem \ref{thm: exist uniq nonlin} is provided by the following lemma:\\
\begin{lemma}\label{lem: Cauchy seq}
For fixed $\mu_0\in\CM^+([0,1])$ and $(\alpha_k)_{k\in\N}$ satisfying Assumption \ref{ass: partition}, the corresponding sequence $(\mu^k)_{k\in\N}$ defined by \eqref{eqn: Euler scheme} is a Cauchy sequence in $C([0,T];\CM^+([0,1]))$. In particular, there is a constant $C$ such that
\begin{equation*}
\sup_{\tau\in[0,T]}\|\mu^k_\tau - \mu^m_\tau\|^*_\BL \leqslant \,C\, \max_{j\in\{0,\ldots,N_k-1\}} \left( t^k_{j+1}-t^k_j \right),
\end{equation*}
for all $k,m\in\N$ satisfying $m\geqslant k$.
\end{lemma}\\

\begin{proof}%  *exactly* the same reasoning works also for partitions \alpha_k and \alpha_m with m>=k! (more natural for Cauchy sequences)
Fix $k,m\in\N$ with $m\geqslant k$, let $\tau\in[0,T]$ be arbitrary and let $j\in\{0,\ldots,N_k-1\}$ be such that $\tau\in(t^k_j,t^k_{j+1}]$. Define, for appropriate $N^{(j)}\geqslant1$, the ordered set
\begin{equation}
\{\tau_\ell\;\,:\;\, 0\leqslant \ell \leqslant N^{(j)}\}:= \{t^k_j\}\, \cup \, \bigg( \alpha_m\cap(t^k_j,t^k_{j+1}] \bigg) \, \cup \, \{t^k_{j+1}\}.\label{eqn: def set tau}
\end{equation}
The set $\alpha_m\cap(t^k_j,t^k_{j+1}]$ contains all $t^m_\ell$, $\ell\in\{1,\ldots,N_m\}$, such that $t^k_j<t^m_\ell\leqslant t^k_{j+1}$. For the sake of being complete, we emphasize that any duplicate elements that might occur on the right-hand side of \eqref{eqn: def set tau} are not `visible' in the set on the left-hand side. Assume that $i\in\{0,\ldots,N^{(j)}-1\}$ is such that $\tau\in(\tau_i,\tau_{i+1}]$. To simplify notation, we write $v^\kappa_\ell:=v[\mu^\kappa_{\tau_\ell}]$ for all $\kappa\in\N$ and $\ell\in\{0,\ldots,N^{(j)}\}$. Define $i_0\in\{0,\ldots,N_m\}$ to be the smallest index such that $t^m_{i_0}\geqslant t^k_j$.

%%%%%%%%%%%%%%%%%%%%%%%%%%%%%%%%%%%%%%%%%%%%%
\paragraph{Case 1: $t^m_{i_0}=t^k_j$}
In this case, there is a $q\in\{0,\ldots,N_m-1\}$ such that $\tau_i=t^m_q$. Hence,
\begin{equation*}
\mu^k_\tau = Q^{v^k_0}_{\tau-\tau_i}\,\mu^k_{\tau_i},\qquad\text{and}\qquad \mu^m_\tau = Q^{v^m_i}_{\tau-\tau_i}\,\mu^m_{\tau_i}.
\end{equation*}
We estimate
\begin{align}
\nonumber \|\mu^k_\tau-\mu^m_\tau\|^*_\BL \leqslant&\, \|Q^{v^k_0}_{\tau-\tau_i}(\mu^k_{\tau_i}-\mu^m_{\tau_i})\|^*_\BL \,+\, \|\big(Q^{v^k_0}_{\tau-\tau_i} - Q^{v^m_i}_{\tau-\tau_i}\big)\mu^m_{\tau_i}\|^*_\BL\\
\nonumber \vspace{0.5cm}\\
\nonumber \leqs&\, \|\mu^k_{\tau_i}-\mu^m_{\tau_i}\|^*_\BL \, \exp\bigg(|v^k_0|_\Lip\,(\tau-\tau_i)+\|f\|_\BL\,(\tau-\tau_i)\,e^{|v^k_0|_\Lip\,(\tau-\tau_i)}\bigg)\\
\nonumber & + \|v^k_0-v^m_i\|_\infty\, \,\|\mu^m_{\tau_i}\|_\TV\,\exp\bigg(L\,(\tau-\tau_i)+\|f\|_\BL\,(\tau-\tau_i)\,e^{L\,(\tau-\tau_i)}\bigg)\cdot\\
& \hspace{3.95cm}\cdot\Big[(\tau-\tau_i)+(\tau-\tau_i)^2\|f\|_\infty \,e^{\|f\|_\infty\,(\tau-\tau_i)}\Big],\label{eqn:estimate muk - muk+1 using Q case1}
\end{align} 
using Corollary \ref{cor: Qmu - Qnu} and Lemma \ref{lem: Qvmu - Qv'mu}. Here, $L$ denotes $\min(|v^k_0|_\Lip\,,|v^m_i|_\Lip)$. In view of Lemma \ref{lem: set timeslices bdd}, we define $R:=\|\mu_0\|_\TV \cdot\exp(\|f\|_\infty\,T)$. From Lemma \ref{lem: Qtmu - Qsmu} (with $s=0$), and Parts (ii) and (iv) %\ref{ass: part Lipsch in x} and \ref{ass: part Lipsch in mu incl R} 
of Assumption \ref{ass: v properties} it follows that
\begin{align}
\nonumber\|v^k_0-v^m_i\|_\infty &\leqslant M_R\,\Big(\|\mu^k_{\tau_0}-\mu^m_{\tau_0}\|^*_\BL + \sum_{\ell=1}^i\|\mu^m_{\tau_\ell}-\mu^m_{\tau_{\ell-1}}\|^*_\BL\Big) \\
\nonumber &\leqslant M_R\,\|\mu^k_{\tau_0}-\mu^m_{\tau_0}\|^*_\BL + M_R\,\sum_{\ell=1}^i\|Q^{v^m_{\ell-1}}_{\tau_\ell-\tau_{\ell-1}}\,\mu^m_{\tau_{\ell-1}}-\mu^m_{\tau_{\ell-1}}\|^*_\BL\\
%\nonumber \leqslant& M_R\,\|\mu^k_{\tau_0}-\mu^m_{\tau_0}\|^*_\BL + M_R\,\sum_{\ell=1}^i\|\mu^{k+1}_{\tau_{\ell-1}}\|_\TV\, \bigl(\|f\|_\infty + \|v^{k+1}_{\ell-1}\|_\infty\bigr)\, e^{\|f\|_\infty\,(\tau_\ell-\tau_{\ell-1})}\,(\tau_\ell-\tau_{\ell-1})\\
\nonumber &\leqslant M_R\,\|\mu^k_{\tau_0}-\mu^m_{\tau_0}\|^*_\BL + M_R \sum_{\ell=1}^i R\,\bigl(\|f\|_\infty + K_R\bigr)\,e^{\|f\|_\infty\,T}\,(\tau_\ell-\tau_{\ell-1})\\
&\leqslant M_R\,\|\mu^k_{\tau_0}-\mu^m_{\tau_0}\|^*_\BL + M_R \,R\,\bigl(\|f\|_\infty + K_R\bigr)\,e^{\|f\|_\infty\,T}\,(\tau_i-\tau_0).\label{eqn:est vk0 - vk+1i infty norm case1}
\end{align}
We combine \eqref{eqn:estimate muk - muk+1 using Q case1} and \eqref{eqn:est vk0 - vk+1i infty norm case1}, and use Part (iii) of Assumption \ref{ass: v properties} and the basic estimates $\tau-\tau_i\leqslant \tau_{i+1}-\tau_i$ and $\tau_{i+1}-\tau_i\leqslant T$ (in suitable places) to obtain that
\begin{align}
\nonumber \|\mu^k_\tau-\mu^m_\tau\|^*_\BL \leqs&\, \exp\big(A_1\,(\tau_{i+1}-\tau_i)\big)\,\|\mu^k_{\tau_i}-\mu^m_{\tau_i}\|^*_\BL \\
\nonumber &+ A_2\,(\tau_{i+1}-\tau_i)\,\|\mu^k_{\tau_0}-\mu^m_{\tau_0}\|^*_\BL \\
&+ A_3\,(\tau_{i+1}-\tau_i)(\tau_i-\tau_0) \label{eqn: est with constants case1}
\end{align}
for some positive constants $A_1$, $A_2$ and $A_3$ that depend on $f$, $T$ and $R$, but not on $i$ or $j$. This upper bound holds for all $\tau\in(\tau_i,\tau_{i+1}]$.

%%%%%%%%%%%%%%%%%%%%%%%%%%%%%%%%%%%%%%%%%%%%%
\paragraph{Case 2: $t^k_j<t^m_{i_0}$ and $i=0$}
Note that $j\neq0$ and $i_0\neq0$ must hold. We recall the notation $v^\kappa_\ell:=v[\mu^\kappa_{\tau_\ell}]$ for all $\kappa\in\N$ and $\ell\in\{0,\ldots,N^{(j)}\}$. In this case,
\begin{equation*}
\mu^k_\tau = Q^{v^k_0}_{\tau-\tau_0}\,\mu^k_{\tau_0},\qquad\text{and}\qquad \mu^m_\tau = Q^{\bar{v}}_{\tau-\tau_0}\,\mu^m_{\tau_0},
\end{equation*}
where $\bar{v}:=v[\mu^m_{t^m_{i_0-1}}]$. Similar to \eqref{eqn:estimate muk - muk+1 using Q case1}, we have
\begin{align}
\nonumber \|\mu^k_\tau-\mu^m_\tau\|^*_\BL \leqslant&\, \|Q^{v^k_0}_{\tau-\tau_0}(\mu^k_{\tau_0}-\mu^m_{\tau_0})\|^*_\BL \,+\, \|\big(Q^{v^k_0}_{\tau-\tau_0} - Q^{\bar{v}}_{\tau-\tau_0}\big)\mu^m_{\tau_0}\|^*_\BL\\
\nonumber \vspace{0.25cm}\\
\nonumber \leqs&\, \|\mu^k_{\tau_0}-\mu^m_{\tau_0}\|^*_\BL \, \exp\bigg(|v^k_0|_\Lip\,(\tau-\tau_0)+\|f\|_\BL\,(\tau-\tau_0)\,e^{|v^k_0|_\Lip\,(\tau-\tau_0)}\bigg)\\
\nonumber & + \|v^k_0-\bar{v}\|_\infty\, \,\|\mu^m_{\tau_0}\|_\TV\,\exp\bigg(L\,(\tau-\tau_0)+\|f\|_\BL\,(\tau-\tau_0)\,e^{L\,(\tau-\tau_0)}\bigg)\cdot\\
& \hspace{3.95cm}\cdot\Big[(\tau-\tau_0)+(\tau-\tau_0)^2\|f\|_\infty \,e^{\|f\|_\infty\,(\tau-\tau_0)}\Big],\label{eqn:estimate muk - muk+1 using Q case2}
\end{align}
where $L=\min(|v^k_0|_\Lip\,,|\bar{v}|_\Lip)$. We define $R:=\|\mu_0\|_\TV \cdot\exp(\|f\|_\infty\,T)$; cf.~Lemma \ref{lem: set timeslices bdd}. %From Lemma \ref{lem: Qtmu - Qsmu} (with $s=0$), Parts (ii) and (iv) %\ref{ass: part Lipsch in x} and \ref{ass: part Lipsch in mu incl R} 
%of Assumption \ref{ass: v properties} it follows that
The analogon of \eqref{eqn:est vk0 - vk+1i infty norm case1} is
\begin{align}
\nonumber\|v^k_0-\bar{v}\|_\infty &\leqslant M_R\,\Big(\|\mu^k_{\tau_0}-\mu^m_{\tau_0}\|^*_\BL + \|\mu^m_{\tau_0}-\mu^m_{\bar{\tau}}\|^*_\BL\Big)\\
\nonumber &= M_R\,\|\mu^k_{\tau_0}-\mu^m_{\tau_0}\|^*_\BL + M_R\,\|Q^{\bar{v}}_{\tau_0-\bar{\tau}}\,\mu^m_{\bar{\tau}}-\mu^m_{\bar{\tau}}\|^*_\BL\\
&\leqslant M_R\,\|\mu^k_{\tau_0}-\mu^m_{\tau_0}\|^*_\BL + M_R \, R\,\bigl(\|f\|_\infty + K_R\bigr)\,e^{\|f\|_\infty\,T}\,(\tau_0-\bar{\tau}),\label{eqn:est vk0 - vk+1i infty norm case2}
\end{align}
with $\bar{\tau}:=t^m_{i_0-1}$. Together \eqref{eqn:estimate muk - muk+1 using Q case2} and \eqref{eqn:est vk0 - vk+1i infty norm case2} yield
\begin{align}
\nonumber \|\mu^k_\tau-\mu^m_\tau\|^*_\BL \leqs&\, \big[\exp\big(A_1\,(\tau_1-\tau_0)\big) + A_2\,(\tau_1-\tau_0) \big]\,\|\mu^k_{\tau_0}-\mu^m_{\tau_0}\|^*_\BL \\
&+ A_3\,(\tau_1-\tau_0)(\tau_0-\bar{\tau})\label{eqn: est with constants case2}
\end{align}
for the same positive constants $A_1$, $A_2$ and $A_3$ as in \eqref{eqn: est with constants case1}. Here, we used Part (iii) of Assumption \ref{ass: v properties} and the estimates $\tau-\tau_0\leqslant \tau_1-\tau_0$ and $\tau_1-\tau_0\leqslant T$. The upper bound \eqref{eqn: est with constants case2} holds for all $\tau\in(\tau_0,\tau_1]$.

%%%%%%%%%%%%%%%%%%%%%%%%%%%%%%%%%%%%%%%%%%%%%
\paragraph{Case 3: $t^k_j<t^m_{i_0}$ and $i\geqslant1$}
In this case, $t^k_j<\tau_i<t^k_{j+1}$ and hence there is a $q\in\{1,\ldots,N_m-1\}$ such that $\tau_i=t^m_q$. We have
\begin{equation*}
\mu^k_\tau = Q^{v^k_0}_{\tau-\tau_i}\,\mu^k_{\tau_i},\qquad\text{and}\qquad \mu^m_\tau = Q^{v^m_i}_{\tau-\tau_i}\,\mu^m_{\tau_i}.
\end{equation*}
Estimate \eqref{eqn:estimate muk - muk+1 using Q case1} also holds in this case. Because $t^m_{i_0}>t^k_j$ there is no $q\in\{0,\ldots,N_m-1\}$ such that $\tau_0=t^m_q$, and therefore $v[\,\cdot\,]$ is not to be evaluated at $\mu^m_{\tau_0}$. Consequently, we have instead of \eqref{eqn:est vk0 - vk+1i infty norm case1},
\begin{align*}
\|v^k_0-v^m_i\|_\infty \leqslant& \|v^k_0-\bar{v}\|_\infty + \|v^m_1-\bar{v}\|_\infty + \sum_{\ell=2}^i\|v^m_{\ell}-v^m_{\ell-1}\|^*_\BL\\
\nonumber \leqslant& \|v^k_0-\bar{v}\|_\infty + M_R\,\|Q^{\bar{v}}_{\tau_1-\bar{\tau}}\,\mu^m_{\bar{\tau}}-\mu^m_{\bar{\tau}}\|^*_\BL  + M_R\,\sum_{\ell=2}^i\|Q^{v^m_{\ell-1}}_{\tau_\ell-\tau_{\ell-1}}\,\mu^m_{\tau_{\ell-1}}-\mu^m_{\tau_{\ell-1}}\|^*_\BL
\end{align*}
with $\bar{v}:=v[\mu^m_{t^m_{i_0-1}}]$ and $\bar{\tau}:=t^m_{i_0-1}$. Note that the sum on the right-hand side might be empty. Using the idea of \eqref{eqn:est vk0 - vk+1i infty norm case1} and the result of \eqref{eqn:est vk0 - vk+1i infty norm case2}, we obtain
\begin{align}
\nonumber \|v^k_0-v^m_i\|_\infty \leqslant&\, M_R\,\|\mu^k_{\tau_0}-\mu^m_{\tau_0}\|^*_\BL + M_R \, R\,\bigl(\|f\|_\infty+ K_R\bigr)\,e^{\|f\|_\infty\,T}\,(\tau_0-\bar{\tau})\\
\nonumber & + M_R \,R\,\bigl(\|f\|_\infty + K_R\bigr)\,e^{\|f\|_\infty\,T}\,(\tau_i-\bar{\tau})\\
\leqslant&\, M_R\,\|\mu^k_{\tau_0}-\mu^m_{\tau_0}\|^*_\BL + 2\,M_R \, R\,\bigl(\|f\|_\infty + K_R\bigr)\,e^{\|f\|_\infty\,T}\,(\tau_i-\bar{\tau}).\label{eqn:est vk0 - vk+1i infty norm case3}
\end{align}
Due to \eqref{eqn:estimate muk - muk+1 using Q case1} and \eqref{eqn:est vk0 - vk+1i infty norm case3}, we have
\begin{align}
\nonumber \|\mu^k_\tau-\mu^m_\tau\|^*_\BL \leqs&\, \exp\big(A_1\,(\tau_{i+1}-\tau_i)\big)\,\|\mu^k_{\tau_i}-\mu^m_{\tau_i}\|^*_\BL \\
\nonumber &+ A_2\,(\tau_{i+1}-\tau_i)\,\|\mu^k_{\tau_0}-\mu^m_{\tau_0}\|^*_\BL \\
&+ 2\, A_3\,(\tau_{i+1}-\tau_i)(\tau_i-\bar{\tau}) \label{eqn: est with constants case3}
\end{align}
for all $\tau\in(\tau_i,\tau_{i+1}]$, where $A_1$, $A_2$ and $A_3$ are the same constants as in \eqref{eqn: est with constants case1} and \eqref{eqn: est with constants case2}.\\
%Therefore, we can take the supremum over $\tau$ on the left-hand side and have the same bound:
%\begin{align*}
%\sup_{\tau\in(\tau_0,\tau_1]}\|\mu^k_\tau-\mu^{k+1}_\tau\|^*_\BL \leqs&\, \big[\exp\big(A_1\,(\tau_1-\tau_0)\big) + A_2\,(\tau_1-\tau_0) \big]\,\|\mu^k_{\tau_0}-\mu^{k+1}_{\tau_0}\|^*_\BL \\
%&+ A_3\,(\tau_1-\tau_0)(\tau_0-\bar{\tau}).
%\end{align*}
\\
%%%%%%%%%%%%%%%%%%%%%%%%%%%%%%%%%%%%%%%%
We now combine the estimates obtained in Cases 1, 2 and 3: it follows from \eqref{eqn: est with constants case1}, \eqref{eqn: est with constants case2} and \eqref{eqn: est with constants case3} that
\begin{align*}
\sup_{\tau\in(\tau_i,\tau_{i+1}]}\|\mu^k_\tau-\mu^m_\tau\|^*_\BL \leqs&\, \exp\big(A_1\,(\tau_{i+1}-\tau_i)\big)\,\sup_{\tau\in(\tau_{i-1},\tau_i]}\|\mu^k_{\tau}-\mu^m_{\tau}\|^*_\BL \\
&+ A_2\,(\tau_{i+1}-\tau_i)\,\|\mu^k_{\tau_0}-\mu^m_{\tau_0}\|^*_\BL \\
&+ 4\,A_3\,M^{(k)}\,(\tau_{i+1}-\tau_i),
\end{align*}
for all $i\in\{1,\ldots,N^{(j)}-1\}$, while for $i=0$
\begin{align*}
\sup_{\tau\in(\tau_0,\tau_1]}\|\mu^k_\tau-\mu^m_\tau\|^*_\BL \leqs&\, \big[\exp\big(A_1\,(\tau_1-\tau_0)\big)+A_2\,(\tau_1-\tau_0)\big]\,\|\mu^k_{\tau_0}-\mu^m_{\tau_0}\|^*_\BL\\
&+ A_3\,M^{(k)}\,(\tau_1-\tau_0).
\end{align*}
We have used that $\tau_i-\tau_0\leqslant M^{(k)}$ in \eqref{eqn: est with constants case1}, $\tau_0-\bar{\tau}\leqslant t^m_{i_0}-t^m_{i_0-1}\leqslant M^{(m)}\leqslant M^{(k)}$ in \eqref{eqn: est with constants case2} and 
\begin{equation*}
\tau_{i}-\bar{\tau}\leqslant\, \tau_{N^{(j)}}-\tau_0 + \tau_0 - \bar{\tau} \leqslant\, M^{(k)} + t^m_{i_0}-t^m_{i_0-1} \leqslant\, M^{(k)} + M^{(m)} \leqslant\, 2\,M^{(k)}
\end{equation*}
in \eqref{eqn: est with constants case3}. This is the place where we use explicitly that partition $\alpha_m$ is `finer' (or: `not coarser') than $\alpha_k$ in the sense that $M^{(m)}\leqslant M^{(k)}$; cf.~Assumption \ref{ass: partition}. By an induction argument one can show that the upper bound
\begin{align}
\nonumber \sup_{\tau\in(\tau_i,\tau_{i+1}]}\|\mu^k_\tau-\mu^m_\tau\|^*_\BL \leqs&\, 
\sum_{\ell=0}^i \left(\prod_{q=\ell+1}^i \exp\big(A_1\,(\tau_{q+1}-\tau_q)\big) \right)\cdot\\
\nonumber & \hspace{2em} \cdot \bigg[A_2\,(\tau_{\ell+1}-\tau_\ell)\,\|\mu^k_{\tau_0}-\mu^m_{\tau_0}\|^*_\BL + 4\,A_3\,M^{(k)}\,(\tau_{\ell+1}-\tau_\ell)\bigg]\\
&+ \, \left(\prod_{q=0}^i \exp\big(A_1\,(\tau_{q+1}-\tau_q)\big) \right)\, \|\mu^k_{\tau_0}-\mu^m_{\tau_0}\|^*_\BL\label{eqn: est after induction within subinterval}
\end{align}
holds for all $i\in\{0,\ldots,N^{(j)}-1\}$. The products in brackets are equal to $\exp\big(A_1\,(\tau_{i+1}-\tau_{\ell+1})\big)$ and $\exp\big(A_1\,(\tau_{i+1}-\tau_{0})\big)$, respectively. By using these explicit expressions and by taking the supremum over $i$ on the left-hand and right-hand sides of \eqref{eqn: est after induction within subinterval}, we obtain
\begin{align}
\nonumber \sup_{\tau\in(\tau_0,\tau_{N^{(j)}}]}\|\mu^k_\tau-\mu^m_\tau\|^*_\BL \leqs&\, 
\sum_{\ell=0}^{N^{(j)}-1} \exp\big(A_1\,(\tau_{N^{(j)}}-\tau_{\ell+1})\big) \cdot\\
\nonumber & \hspace{2em} \cdot \big[A_2\,(\tau_{\ell+1}-\tau_\ell)\,\|\mu^k_{\tau_0}-\mu^m_{\tau_0}\|^*_\BL + 4\,A_3\,M^{(k)}\,(\tau_{\ell+1}-\tau_\ell)\big]\\
&+ \, \exp\big(A_1\,(\tau_{N^{(j)}}-\tau_0)\big) \, \|\mu^k_{\tau_0}-\mu^m_{\tau_0}\|^*_\BL.\label{eqn: est after sup over i and product within subinterval}
\end{align}
Since $\tau_{N^{(j)}}-\tau_{\ell+1}\leqslant \tau_{N^{(j)}}-\tau_0$ for all $\ell\in\{0,\ldots,N^{(j)}-1\}$, and $\tau_0=t^k_j$ and $\tau_{N^{(j)}}=t^k_{j+1}$, it follows from \eqref{eqn: est after sup over i and product within subinterval} that
\begin{align}
\nonumber \sup_{\tau\in(t^k_j,t^k_{j+1}]}\|\mu^k_\tau-\mu^m_\tau\|^*_\BL \leqs&\, 
\bigg[1+  A_2\,\sum_{\ell=0}^{N^{(j)}-1} (\tau_{\ell+1}-\tau_\ell)\bigg]\cdot \exp\big(A_1\,(t^k_{j+1}-t^k_j)\big) \cdot \|\mu^k_{t^k_j}-\mu^m_{t^k_j}\|^*_\BL\\
\nonumber & + 4\,A_3\,M^{(k)}\,(t^k_{j+1}-t^k_j)\,\exp\big(A_1\,(t^k_{j+1}-t^k_j)\big)\, \sum_{\ell=0}^{N^{(j)}-1}(\tau_{\ell+1}-\tau_\ell)\\
\nonumber =&\, 
\big[1+  A_2\,(t^k_{j+1}-t^k_j)\big]\cdot\exp\big(A_1\,(t^k_{j+1}-t^k_j)\big) \cdot \|\mu^k_{t^k_j}-\mu^m_{t^k_j}\|^*_\BL\\
 & + 4\,A_3\,M^{(k)}\,(t^k_{j+1}-t^k_j)\,\exp\big(A_1\,(t^k_{j+1}-t^k_j)\big).
\end{align}
Hence, we have that
\begin{align}
\nonumber \sup_{\tau\in(t^k_j,t^k_{j+1}]}\|\mu^k_\tau-\mu^m_\tau\|^*_\BL \leqs&\, 
\big[1+  A_2\,(t^k_{j+1}-t^k_j)\big]\cdot\exp\big(A_1\,(t^k_{j+1}-t^k_j)\big) \cdot \sup_{\tau\in(t^k_{j-1},t^k_j]}\|\mu^k_\tau-\mu^m_\tau\|^*_\BL\\
 & + 4\,A_3\,M^{(k)}\,(t^k_{j+1}-t^k_j)\,\exp\big(A_1\,(t^k_{j+1}-t^k_j)\big)
\end{align}
for all $j\in\{1,\ldots,N_k-1\}$, and for $j=0$ we have
\begin{equation}
\sup_{\tau\in(t^k_0,t^k_1]}\|\mu^k_\tau-\mu^m_\tau\|^*_\BL \leqs\, 4\,A_3\,M^{(k)}\,(t^k_1-t^k_0)\,\exp\big(A_1\,(t^k_1-t^k_0)\big),
\end{equation}
because $\mu^k_{t^k_0}=\mu_0=\mu^m_{t^k_0}$. By an induction argument similar to the one leading to \eqref{eqn: est after induction within subinterval}, we obtain that
\begin{align}
\nonumber \sup_{\tau\in(t^k_j,t^k_{j+1}]}\|\mu^k_\tau-\mu^m_\tau\|^*_\BL \leqs&\, 
\sum_{\ell=0}^j \left(\prod_{q=\ell+1}^j \big[1+  A_2\,(t^k_{q+1}-t^k_q)\big]\cdot\exp\big(A_1\,(t^k_{q+1}-t^k_q)\big) \right)\cdot\\
& \hspace{4em} \cdot 4\,A_3\,M^{(k)}\,(t^k_{\ell+1}-t^k_\ell)\,\exp\big(A_1\,(t^k_{\ell+1}-t^k_\ell)\big)\label{eqn: est after induction interval j}
\end{align}
for all $j\in\{0,\ldots,N_k-1\}$. Note that $\big[1+  A_2\,(t^k_{q+1}-t^k_q)\big]\leqslant\exp\big(A_2\,(t^k_{q+1}-t^k_q)\big)$ for all $q\in\{0,\ldots,N_k-1\}$. Define $A_4:=A_1+A_2$. It follows from \eqref{eqn: est after induction interval j} that
\begin{align}
\nonumber \sup_{\tau\in(t^k_j,t^k_{j+1}]}\|\mu^k_\tau-\mu^m_\tau\|^*_\BL \leqs&\, 
4\,A_3\,M^{(k)}\,\exp(A_1\,T)\,\sum_{\ell=0}^j \left(\prod_{q=\ell+1}^j \exp\big(A_4\,(t^k_{q+1}-t^k_q)\big)\right) \cdot (t^k_{\ell+1}-t^k_\ell)\\
\nonumber \leqs&\, 
4\,A_3\,M^{(k)}\,\exp(A_1\,T)\,\sum_{\ell=0}^j \exp\big(A_4\,(t^k_{j+1}-t^k_{\ell+1})\big) \cdot (t^k_{\ell+1}-t^k_\ell)\\
\leqs&\, 
4\,A_3\,M^{(k)}\,\exp((A_1+A_4)\,T)\,(t^k_{j+1}-t^k_0).\label{eqn: further est after induction interval j}
\end{align}
We take the supremum over $j$ on both sides of the inequality \eqref{eqn: further est after induction interval j} and get
\begin{equation}
\sup_{\tau\in[0,T]}\|\mu^k_\tau-\mu^m_\tau\|^*_\BL \leqs\, 
4\,A_3\,M^{(k)}\,T\,\exp((A_1+A_4)\,T).\label{eqn: final est Cauchy seq}
\end{equation}
Note that we extended the supremum from $\tau\in(0,T]$ to $\tau\in[0,T]$, but this does not change the upper bound. Define $C:=4\,A_3\,T\,\exp((A_1+A_4)\,T)$ to get the result of the lemma. Because $M^{(k)}\to 0$ as $k\to \infty$, the estimate \eqref{eqn: final est Cauchy seq} implies that $(\mu^k)_{k\in\N}$ is a Cauchy sequence.
\end{proof}\\

\begin{remark}
It is crucial that we use $\|\cdot\|_\BL^*$ and not $\|\cdot\|_\TV$ in Lemma \ref{lem: Cauchy seq}. The factor $\|v^k_0-v^m_i\|_\infty$ appears in \eqref{eqn:estimate muk - muk+1 using Q case1} due to Lemma \ref{lem: Qvmu - Qv'mu}. Due to Assumption \ref{ass: v properties}(iv) %\ref{ass: part Lipsch in mu incl R} 
we subsequently obtain an estimate in which $\|\mu^k_{\tau_0}-\mu^m_{\tau_0}\|^*_\BL$ appears. % \textit{as a factor}. This is essential for the result of the lemma. 
Analogous estimates apply to $\|v^k_0-\bar{v}\|_\infty$ in \eqref{eqn:estimate muk - muk+1 using Q case2}. Note that Lemma \ref{lem: Qvmu - Qv'mu} builds on Lemma \ref{lem: Pv mu - Pv' mu}. In \eqref{eqn: pair Pvt-Pv't phi Lipschitz} the Lipschitz property of the test functions is explicitly used and hence, there is no direct way to formulate the result of Lemma \ref{lem: Pv mu - Pv' mu} in terms of $\|\cdot\|_\TV$. Consequently, we do not have an estimate of $\|\mu_\tau^k-\mu_\tau^{k+1}\|_\TV$ against $\|v-v'\|_\infty$ comparable to \eqref{eqn:estimate muk - muk+1 using Q case1}.
\end{remark}\\
\\
%We recall and prove Theorem \ref{thm: exist uniq nonlin}.
%\begin{reptheorem}{thm: exist uniq nonlin}
%For each $\mu_0\in\CM^+([0,1])$ and $\map{v}{\CM([0,1])\times[0,1]}{\R}$ satisfying Assumption \ref{ass: v properties}, there is a unique element of $C([0,T]; \CM^+([0,1]))$ with initial condition $\mu_0$, that is a mild solution\index{mild solution} in the sense of Definition \ref{def: mild soln nonlin}.
%\end{reptheorem}
We are now ready to prove Theorem \ref{thm: exist uniq nonlin}.
\begin{proof}
By definition, $\CMc([0,1])$ is complete in the metric induced by the norm $\|\cdot\|^*_\BL$. %In Lemma \ref{lem: pos meas closed}, 
The space $\CM^+([0,1])$ is a closed subspace of $\CMc([0,1])$, %\\ \textbf{\color{red}(To be checked! I can prove this for TV-norm but not for BL-norm).}\\
so $\CM^+([0,1])$ is complete. %\textbf{\color{red}In \cite{Gwiazda1} this is mentioned underneath Def 2.5, also mentioned in Sander's work.}
%; cf.~Lemma \ref{lem:closedSubsetComplete} in Appendix \ref{app:completeness}. 
Hence, %due to Theorem \ref{thm:completenessAppendix}, 
the space
\[\{\nu\in C([0,T]; \CM^+([0,1]))\,:\, \nu(0)=\mu_0\}\]
is complete for the metric defined for all $\mu,\nu\in C([0,T]; \CM^+([0,1]))$ by \eqref{eqn: metric sup BL norm}.\\
For each initial measure $\mu_0\in\CM^+([0,1])$ and for each $k\in\N$, consider the Euler approximation $\mu^k$ defined by \eqref{eqn: Euler scheme} corresponding to partition $\alpha_k$. This approximation $\mu^k$ is an element of $C([0,T];\CM^+([0,1]))$, because the semigroup $(Q^v_t)_{t\geqs0}$ preserves positivity for all $v\in\BL([0,1])$; see \cite[Corollary 3.4]{EvHiMu_JDE}. In Lemma \ref{lem: Cauchy seq}, we showed that for given $(\alpha_k)_{k\in\N}$ the sequence $(\mu^k)_{k\in\N}$ is a Cauchy sequence in $\{\nu\in C([0,T]; \CM^+([0,1]))\,:\, \nu(0)=\mu_0\}$, which is a complete space, as was argued above. Hence, the sequence $(\mu^k)_{k\in\N}$ converges in $\{\nu\in C([0,T]; \CM^+([0,1]))\,:\, \nu(0)=\mu_0\}$.\\
% and we have existence of a mild solution. Each subsequence as mentioned in Definition \ref{def: mild soln nonlin} converges to that same limit, thus uniqueness of the mild solution holds.
\\
The limit is independent of the sequence of partitions chosen from the class characterized by Assumption \ref{ass: partition}. If $(\alpha_k)_{k\in\N}$ and $(\beta_k)_{k\in\N}$ are two such sequences, then it is possible to construct a sequence $(\gamma_k)_{k\in\N}$ that has a subsequence that is also a subsequence of $(\alpha_k)_{k\in\N}$, and that has (another) subsequence that is a subsequence of $(\beta_k)_{k\in\N}$. Moreover, $(\gamma_k)_{k\in\N}$ can be constructed such that the corresponding sequence of maximal interval lengths is nondecreasing.\\
Let $(\mu^{\alpha_k})_{k\in\N}$, $(\mu^{\beta_k})_{k\in\N}$ and $(\mu^{\gamma_k})_{k\in\N}$ be the corresponding sequences of Euler approximations. The sequence $(\mu^{\gamma_k})_{k\in\N}$ can be shown to converge to the same limit as $(\mu^{\alpha_k})_{k\in\N}$, and to the same limit as $(\mu^{\beta_k})_{k\in\N}$. Hence, $(\mu^{\alpha_k})_{k\in\N}$ and $(\mu^{\beta_k})_{k\in\N}$ converge to the same limit. This finishes the proof.
\end{proof}\\
\\
The proof of Corollary \ref{cor: global exist uniq} builds on the result of Theorem \ref{thm: exist uniq nonlin}.
\begin{proof}
Fix $t\geqslant0$ and let $T>0$ be such that $t\in[0,T]$. For given $\mu_0\in\CM^+([0,1])$ and $\map{v}{\CM([0,1])\times[0,1]}{\R}$ satisfying Assumption \ref{ass: v properties}, a unique mild solution $\mu\in C([0,T];\CM^+([0,1]))$ exists, hence $\mu_t$, the solution at time $t$, exists. We now show that this $\mu_t$ is independent of the choice of $T$.\\
Let $T_1,T_2>0$ and assume without loss of generality that $T_1<T_2$. For the given $\mu_0\in\CM^+([0,1])$ and $\map{v}{\CM([0,1])\times[0,1]}{\R}$, let $\mu$ denote the mild solution in $C([0,T_1];\CM^+([0,1]))$ obtained by partitioning $[0,T_1]$. Take a sequence of partitions $(\alpha_k)_{k\in\N}\subset[0,T_1]$ satisfying Assumption \ref{ass: partition}, with corresponding Euler approximations $(\mu^k)_{k\in\N}$. Next, construct a sequence of partitions $(\beta_k)_{k\in\N}\subset[0,T_2]$ satisfying Assumption \ref{ass: partition}, such that $\alpha_k\subset\beta_k$ for each $k\in\N$. More specifically, restricted to $[0,T_1]$ each partition $\beta_k$ coincides with $\alpha_k$. Note that such $(\beta_k)_{k\in\N}$ exists. Let $(\nu^k)_{k\in\N}$ be the sequence of Euler approximations corresponding to $(\beta_k)_{k\in\N}$.\\
For each $k\in\N$, the restriction $\nu^k\big|_{[0,T_1]}$ is defined by \eqref{eqn: Euler scheme} with respect to the partition $\beta_k \cap [0,T_1] = \alpha_k$. Hence $\nu^k\big|_{[0,T_1]}$ is defined in the same way as $\mu^k$, and thus $\sup_{\tau\in[0,T_1]}\|\mu^k_\tau-\nu^k_\tau\|^*_\BL=0$ or simply $\nu^k\big|_{[0,T_1]}=\mu^k$. Consequently, the same must hold in the limit as $k\to\infty$, because of the triangle inequality: 
\begin{align*}
\sup_{\tau\in[0,T_1]}\|\mu_\tau-\nu_\tau\|^*_\BL \leqs&\, \sup_{\tau\in[0,T_1]}\|\mu^k_\tau-\nu^k_\tau\|^*_\BL + \sup_{\tau\in[0,T_1]}\|\mu^k_\tau-\mu_\tau\|^*_\BL + \sup_{\tau\in[0,T_1]}\|\nu^k_\tau-\nu_\tau\|^*_\BL\\
\leqs&\, \underbrace{\sup_{\tau\in[0,T_1]}\|\mu^k_\tau-\nu^k_\tau\|^*_\BL}_{=0} + \underbrace{\sup_{\tau\in[0,T_1]}\|\mu^k_\tau-\mu_\tau\|^*_\BL}_{\to0} + \underbrace{\sup_{\tau\in[0,T_2]}\|\nu^k_\tau-\nu_\tau\|^*_\BL}_{\to0}.
\end{align*}
So, $\sup_{\tau\in[0,T_1]}\|\mu_\tau-\nu_\tau\|^*_\BL=0$. Hence, $\mu_\tau=\nu_\tau$ for all $\tau\in[0,T_1]$ and thus the solution at time $\tau$ is independent of the final time chosen. 
\end{proof}\\
\\% for cont dep, don't need the factor 4
Finally, we prove Theorem \ref{thm: cont dep}:
\begin{proof}
Let the mild solutions $\mu$ and $\nu$ be given and let $(\alpha_k)_{k\in\N}$ be an arbitrary sequence of partitions of $[0,T]$ satisfying Assumption \ref{ass: partition}. Let $(\mu^k)_{k\in\N}$ and $(\nu^k)_{k\in\N}$ denote the sequences of Euler approximations defined by \eqref{eqn: Euler scheme}, both for the sequence of partitions $(\alpha_k)_{k\in\N}$, and with initial conditions $\mu_0$ and $\nu_0$, respectively.\\
Since $\mu$ and $\nu$ are mild solutions
\begin{align*}
\mu = \lim_{k\to\infty} \mu^k,&\qquad \text{and} \\
\nu = \lim_{k\to\infty} \nu^k&
\end{align*}
hold, with convergence in the metric \eqref{eqn: metric sup BL norm}. It follows from Lemma \ref{lem: set timeslices bdd} that all elements of
\begin{equation*}
\{ \mu^k_t: k\in\N, t\in[0,T] \}\cup\{ \nu^k_t: k\in\N, t\in[0,T] \}
\end{equation*}
are bounded by $R:=\tilde{R}\,\exp(\|f\|_\infty\,T)$ in both $\|\cdot\|_\TV$ and $\|\cdot\|^*_\BL$. Fix $k\in\N$, let $\alpha_k:=\{t^k_0,\ldots,t^k_{N_k}\}$ and take $j$ such that $\tau\in(t^k_j,t^k_{j+1}]$. Consider the difference $\|\mu^k_\tau-\nu^k_\tau\|^*_\BL$.\\ 
\\
We use an estimate in the spirit of \eqref{eqn:estimate muk - muk+1 using Q case1}--\eqref{eqn:est vk0 - vk+1i infty norm case1}--\eqref{eqn: est with constants case1}. Note that the proof of Lemma \ref{lem: Cauchy seq} also holds if $k=m$, which implies $N^{(j)}=1$ and hence $i=0$. It follows from \eqref{eqn:estimate muk - muk+1 using Q case1}--\eqref{eqn:est vk0 - vk+1i infty norm case1}, with $i=0$ and with $\nu^k$ instead of $\mu^m$, that
\begin{equation}
\|\mu^k_\tau-\nu^k_\tau\|^*_\BL \leqs\, \big[1 +  B_2\,(t^k_{j+1}-t^k_j) \big]\,\exp\bigg(B_1\,(t^k_{j+1}-t^k_j)\bigg)\,\|\mu^k_{t^k_j}-\nu^k_{t^k_j}\|^*_\BL \label{eqn: est cont dep nonlin before recurs}
\end{equation}
for some positive constants $B_1$ and $B_2$ that depend on $f$, $T$ and $\tilde{R}$, but not on $j$ or $k$. This estimate holds for all $\tau\in(t^k_j,t^k_{j+1}]$ and resembles \eqref{eqn: est with constants case1}. We take the supremum over $\tau\in(t^k_j,t^k_{j+1}]$ on the left-hand side of \eqref{eqn: est cont dep nonlin before recurs}, apply this relation recursively and take the supremum over $j$ to obtain that
\begin{align}
\nonumber \sup_{\tau\in[0,T]}\|\mu^k_\tau-\nu^k_\tau\|^*_\BL \leqslant&  \left(\prod_{\ell=0}^{N_k-1} \big[1 +  B_2\,(t^k_{\ell+1}-t^k_\ell) \big]\,\exp\big(B_1\,(t^k_{\ell+1}-t^k_\ell)\big)\right)\,\|\mu_{0}-\nu_{0}\|^*_\BL\\
\nonumber \leqs& \left(\prod_{\ell=0}^{N_k-1} \exp\big(B_2\,(t^k_{\ell+1}-t^k_\ell) \big)\,\exp\big(B_1\,(t^k_{\ell+1}-t^k_\ell)\big)\right)\,\|\mu_{0}-\nu_{0}\|^*_\BL\\
\nonumber\leqs&\, \exp\big((B_1+B_2)\,(t^k_{N_k}-t^k_0) \big)\,\|\mu_{0}-\nu_{0}\|^*_\BL\\
=&\, \exp((B_1+B_2)\,T )\,\|\mu_{0}-\nu_{0}\|^*_\BL,\label{eqn: est cont dep nonlin after recurs}
\end{align}
for all $k\in\N$. The triangle inequality yields
\begin{equation*}
\sup_{\tau\in[0,T]}\|\mu_\tau-\nu_\tau\|^*_\BL \leqs \sup_{\tau\in[0,T]}\|\mu^k_\tau-\nu^k_\tau\|^*_\BL + \underbrace{\sup_{\tau\in[0,T]}\|\mu^k_\tau-\mu_\tau\|^*_\BL}_{\to0} + \underbrace{\sup_{\tau\in[0,T]}\|\nu^k_\tau-\nu_\tau\|^*_\BL}_{\to0},
\end{equation*}
whence the same estimate as in \eqref{eqn: est cont dep nonlin after recurs} holds for $\sup_{\tau\in[0,T]}\|\mu_\tau-\nu_\tau\|^*_\BL$.
\end{proof}\\

\begin{remark}
We would have been inclined to use directly \eqref{eq:cont dep initial cond} on the interval $(t^k_j,t^k_{j+1}]$, %the result of Proposition \ref{prop:cont dependence}, 
instead of deriving \eqref{eqn: est cont dep nonlin before recurs}. We need, however, the exact dependence on $(t^k_{j+1}-t^k_j)$ of the prefactor, to make sure that -- after iteration over $j$ -- the prefactor in \eqref{eqn: est cont dep nonlin after recurs} is bounded. This dependence is not (directly) provided by \eqref{eq:cont dep initial cond}, nor by the proof of \cite[Proposition 3.5]{EvHiMu_JDE}.%Proposition \ref{prop:cont dependence}.
\end{remark}\\

\begin{remark}\label{rem:ch 5 Gronwall is used}
The result of Theorem \ref{thm: cont dep} relies -- via Corollary \ref{cor: Qmu - Qnu} and Lemma \ref{lem: Qvmu - Qv'mu} -- on Gronwall's Inequality.
%We used an argument involving Gronwall's Inequality\index{Gronwall's Lemma} to obtain the result of Theorem \ref{thm: cont dep}.
This is possible here because we restrict ourselves to Lipschitz perturbations. In our previous work \cite{EvHiMu_JDE} we considered the more general class of \textit{piecewise} bounded Lipschitz perturbations. Hence, there we stated explicitly (see the paragraph before \cite[Proposition 3.5]{EvHiMu_JDE}) that the standard approach did not work.
\end{remark}

%\section{Numerical illustration}
%\textbf{\color{red}Numerics to be done; phase plane}

\section{Discussion}\label{sec: discussion bc meas-dep v}
%SAID IN CH 4:
%The system is not mass-preserving, and the total mass might increase (without clear \textit{a priori} bound).
%This is also why the BL norm is used and not the Wasserstein distance (well-established/holy). Michiel keeps track of the lost mass. But this does not allow for net influx.
%Mild -> weak, lukte niet? term die niet verdwijnt
%{\color{red} Maybe the title is not optimal, since we cannot derive real boundary conditions, only gating in bdry layer}

In this paper we have generalized the results of \cite{EvHiMu_JDE} to measure-dependent velocity fields via a forward-Euler-like approach. Our motivation was to derive flux boundary conditions for situations in which the dynamics are driven by interactions. Such dynamics are in general more interesting than the dynamics that follow from prescribed velocity fields as in \cite{EvHiMu_JDE}. We managed to obtain a converging procedure, but only for bounded Lipschitz continuous right-hand sides. Hence, compared to \cite{EvHiMu_JDE}, our results hold e.g.~for boundary layers in which mass decays, but not for the limit case of vanishing boundary layer. We start off this discussion section (see \S \ref{sec: disc f}) by commenting on the possibility to extend to piecewise bounded Lipschitz right-hand sides and to obtain the limit of vanishing boundary layer. Secondly, we point out (in \S \ref{sec: disc partitions}) how this paper generalizes the results of \cite[Chapter 5]{Evers_PhD} and how a number of open problems mentioned in \cite[Section 5.5]{Evers_PhD} are now resolved. Ultimately, we suggest possible future research (\S \ref{sec: disc future}).

\subsection{Piecewise bounded Lipschitz perturbations}\label{sec: disc f}
%In the proofs of Lemma \ref{lem: Cauchy seq} and Theorem \ref{thm: cont dep} we explicitly used the assumption that the perturbation $f$ is bounded Lipschitz on $[0,1]$. We would have liked to obtain these results for piecewise bounded Lipschitz $f$, in particular to model mass decay at one of the boundaries only. In Chapter \ref{ch:bc pres v} we circumvent the arising problems by providing the solution explicitly. In the current setting, this explicit form would be given for each interval $(t^k_j,t^k_{j+1}]$, $k\in\N$, in \eqref{eqn: Euler scheme} by
To obtain the technical results in \S \ref{sec: summ tech Q}, we explicitly used the assumption that the perturbation $f$ is bounded Lipschitz on $[0,1]$. Theorem \ref{thm: exist uniq nonlin} and Theorem \ref{thm: cont dep} rely on the results in \S \ref{sec: summ tech Q}. We would have liked to obtain these results for piecewise bounded Lipschitz $f$, in particular to model decay of mass at one of the boundaries only (cf.~\cite{EvHiMu_JDE}). In \cite{EvHiMu_JDE} we circumvent the arising problems by providing the solution explicitly in \cite[Proposition 3.3]{EvHiMu_JDE}. In the setting of the present paper, this explicit form would be given for each interval $(t^k_j,t^k_{j+1}]$, $k\in\N$, in \eqref{eqn: Euler scheme} by
\begin{equation}\label{eqn: expl soln nonlin}
\mu^k_t := \Int{[0,1]}{}{\exp\left(\Int{0}{t-t^k_j}{f(\Phi^{v^k_j}_s(x))}{ds}\right)\cdot \delta_{\Phi^{v^k_j}_{t-t^k_j}(x)}}{d\mu^k_{t^k_j}(x)},
\end{equation}
where $v^k_j:=v[\mu^k_{t^k_j}]$. In \cite{EvHiMu_JDE} we showed that it is possible to obtain the estimates needed to establish continuous dependence on initial data, because this explicit form has a regularizing effect on $f$ and its discontinuities due to the integration in time. The key ingredient there, which is absent in the approach of the present work, is the fact that the velocity field is the same for all time. If one wants to prove Theorem \ref{thm: cont dep} using \eqref{eqn: expl soln nonlin} instead of the properties of the semigroup $Q$, one encounters that at some point for any $\Delta t>0$ fixed a Lipschitz estimate of the form
\begin{equation}\label{eqn: Lipsch est int f nonlin}
\|\Int{0}{\Delta t}{f(\Phi^{v}_s(\cdot))}{ds} - \Int{0}{\Delta t}{f(\Phi^{u}_s(\cdot))}{ds}\|_\infty \leqslant C\, \|u-v\|_\infty
\end{equation}
is required, for all $u$ and $v$ taken from a class of admissible %measure-dependent
velocity fields. One would then proceed to estimate $\|u-v\|_\infty$ against the $\BL$-distance of the corresponding measures, using Part (iv) of Assumption \ref{ass: v properties}.\\
\\
In view of \cite{EvHiMu_JDE}, the restriction that the velocity should not be zero at discontinuities of $f$ is reasonable, but even if we are willing to obey that condition, an estimate like \eqref{eqn: Lipsch est int f nonlin} cannot be expected to hold. Let $f(x)=0$ if $x\in[0,1)$ and $f(1)=-1$. Take $\eps>0$ and take $v\equiv\eps$, $u\equiv-\eps$. Then (for $\eps<1/\Delta t$)
\begin{multline*}
\|\Int{0}{\Delta t}{f(\Phi^{v}_s(\cdot))}{ds} - \Int{0}{\Delta t}{f(\Phi^{u}_s(\cdot))}{ds}\|_\infty \geqs |\Int{0}{\Delta t}{f(\Phi^{v}_s(1))}{ds} - \Int{0}{\Delta t}{f(\Phi^{u}_s(1))}{ds}|\\
 = |\Int{0}{\Delta t}{f(1)}{ds} - \Int{0}{\Delta t}{f(1-\eps\,s)}{ds}| = \Delta t.
\end{multline*}
Since $\Delta t>0$ is fixed and $\|u-v\|_\infty=2\eps$ can be made arbitrarily small, \eqref{eqn: Lipsch est int f nonlin} cannot be satisfied.\\
An additional difficulty is that it remains to be seen how we can assure that a condition like $v(1)\neq0$ is satisfied by a velocity field that depends on the solution itself. %{\color{red}Cf. the remark in the intro of this chapter on Ron's condition of inward pointing velocity.}

\subsection{Uniqueness of mild solutions and generality of partitions}\label{sec: disc partitions}
In \cite[Section 5.5]{Evers_PhD} we point out that there are two reasons why we obtained uniqueness of mild solutions there. On the one hand, this is because the constructed approximating sequence converges, thus inevitably each subsequence (cf.~Definition \ref{def: mild soln nonlin}) converges to the same limit. This statement still holds true for the present work. On the other hand, uniqueness holds in \cite[Chapter 5]{Evers_PhD} because there we only constructed one approximating sequence, namely by partitioning the interval $[0,T]$ into $2^k$ subintervals. In this respect, the present paper is a considerable improvement. The class of admissible partitions (see Assumption \ref{ass: partition}) includes partitions into $q^k$ equal subintervals for arbitrary $q\in\Np$; see Example \ref{ex: partitions}. We conjectured in \cite[Section 5.5]{Evers_PhD} that the sequence of corresponding Euler approximations converges, and the results of this paper confirm that conjecture. The fact that, in this case, each interval $(t^k_j,t^k_{j+1}]$ is split into $q^{m-k}$ subintervals $(t^m_\ell,t^m_{\ell+1}]$ is generically treated by introducing the number $N^{(j)}$ and using a recursion over index $i\in\{0,\ldots,N^{(j)}-1\}$ to obtain \eqref{eqn: est after sup over i and product within subinterval}. In \cite[Chapter 5]{Evers_PhD}, however, we performed explicit calculations, using that each $(t^k_j,t^k_{j+1}]$ is split into two subintervals.\\
In \cite[Section 5.5]{Evers_PhD} anticipated that using a sequence of non-uniform partitions of $[0,T]$, would imply the need for a condition regularizing the variation in subinterval lengths to make sure that \textit{all} subintervals become small sufficiently fast as $k\to\infty$. In the present work we show that it suffices to have for the maximum subinterval length $M^{(k)}\to 0$ as $k\to\infty$.\\
The iterative argument in \cite[Chapter 5]{Evers_PhD} requires that the partition for index $k+1$ is a refinement of the partition for index $k$ (more particularly: a division of each subinterval into two). The complications expected to occur if subsequent partitions are not refinements are resolved in the current work, by introduction of the index $i_0$ in the proof of Lemma \ref{lem: Cauchy seq} and allowing for the case $t^k_j\neq t^m_{i_0}$.\\
The final contribution of the present work to be mentioned here is that in Theorem \ref{thm: exist uniq nonlin} we have positively answered the question posed in \cite[Section 5.5]{Evers_PhD} whether the mild solutions obtained as limits of distinct sequences of partitions are actually identical. 

\subsection{Future directions}\label{sec: disc future}
The extension of the results stated in \S \ref{sec: gen meas-dep v} to functions $f$ with discontinuities would clear the way for an approximation procedure like the one treated in \cite{EvHiMu_JDE}. That is, to have $f$ nonzero \emph{only} on the boundary of the domain and to approximate it with a sequence of bounded Lipschitz functions $(f_n)_{n\in\N}\subset\BL([0,1])$. In \cite{EvHiMu_JDE} we showed convergence of the corresponding solutions as $n\to\infty$ (for $v\in\BL([0,1])$ fixed). The challenge would be (i): to establish the well-posedness of the problem for discontinuous $f$, and (ii): to show that the Euler approximation limit and the boundary layer limit commute.\\
\\
Let us focus on the vanishing boundary layer like in \cite{EvHiMu_JDE}. Assume there are regions around $0$ and $1$ in which mass decays, and that these regions shrink to zero width. That is, there is a sequence $(f_n)_{n\in\N}\subset\BL([0,1])$ and there is an $f$ satisfying $f(x)=0$ if $x\in(0,1)$ and e.g.~$f(0)=f(1)=-1$, such that $f_n\to f$ pointwise, and the Lebesgue measure of the set $\{x\in[0,1]: f_n(x)\neq f(x)\}$ tends to zero as $n\to\infty$. If we assume that we can extend the results of this paper to piecewise bounded Lipschitz $f$, then well-posedness for the limit case is guaranteed. It remains to be proven however that the solution for finite boundary layer actually converges to the solution of the limit problem.\\
This is the same question as asking whether the two limits that we take, actually commute. The first limit is in the forward-Euler-like approach to obtain a mild solution. We assigned an index $k$ to the elements in the approximating sequence and proved in Theorem \ref{thm: exist uniq nonlin} that the limit ``$\lim_{k\to\infty}$'' exists (for $f\in\BL([0,1])$). The second limit ``$\lim_{n\to\infty}$'' is the one  involving the sequence $(f_n)_{n\in\N}\subset\BL([0,1])$. Proving the well-posedness for $f$ piecewise bounded Lipschitz, is the same as proving that the limit ``$\lim_{k\to\infty}\lim_{n\to\infty}$'' exists. Proving that the sequence of solutions corresponding to each $f_n$ actually converges to some limit in $C([0,T];\CM^+([0,1]))$ is equivalent to proving that ``$\lim_{n\to\infty}\lim_{k\to\infty}$'' exists. To conclude that the two limits commute, an additional argument is needed. It requires a characterization of ``$\lim_{n\to\infty}\lim_{k\to\infty}$'' that can be compared to ``$\lim_{k\to\infty}\lim_{n\to\infty}$''. Both proving that ``$\lim_{n\to\infty}\lim_{k\to\infty}$'' exists and characterizing the limit can be a difficult task, however, since our current results do not provide an explicit expression for ``$\lim_{k\to\infty}$''. A possible way to characterize the limit ``$\lim_{k\to\infty}$'' could be to show that the mild solution obtained in this paper is actually a weak solution, and use the weak formulation of \eqref{eqn:main equation nonlin} as a characterization. If the solutions obtained in this paper are weak solutions, this is also a further justification of the terminology `mild solutions'.\\
\\
An additional result to be derived concerns the stability with respect to parameters, in particular with respect to $f$ and the specific form of $v$. Stability statements are essential in view of parameter identification. It is important to know how measurement errors in the parameters affect the solution of our model. In fact, Lemma \ref{lem: Qvmu - Qv'mu} already provides stability in $v$ for the solution of \cite{EvHiMu_JDE}, provided that $f\in\BL([0,1])$.\\
Moreover, we would like to study the long-term dynamics of the solutions $t\mapsto \mu_t$ for various initial conditions.

%One step beyond the measure-dependent velocity fields proposed in this paper, is to add anisotropy due to a field of vision. This would link the approach of Chapters \ref{ch:bc pres v} and \ref{ch:bc meas-dep v} to Chapter \ref{ch: aniso v}.\\
%\\

\section*{Acknowledgements}
Until 2015 J.H.M. Evers was a member of the Centre for Analysis, Scientific computing and Applications, and the Institute for Complex Molecular Systems (ICMS) at Eindhoven University of Technology, supported by the Netherlands Organisation for Scientific Research (NWO), Graduate Programme 2010.
%\ldots

\bibliographystyle{alpha}%siam.bst
\bibliography{refs}

\newcommand{\etalchar}[1]{$^{#1}$}
\begin{thebibliography}{GLMC10}

\bibitem[AI05]{Ackleh1}
A.S. Ackleh and K.~Ito.
\newblock Measure-valued solutions for a hierarchically size-structured
  population.
\newblock {\em Journal of Differential Equations}, 217(2):431 -- 455, 2005.

\bibitem[BGCG06]{BenzoniColomboGwiazda}
S.~Benzoni-Gavage, R.M. Colombo, and P.~Gwiazda.
\newblock Measure valued solutions to conservation laws motivated by traffic
  modelling.
\newblock {\em Proceedings of the Royal Society of London A: Mathematical,
  Physical and Engineering Sciences}, 462(2070):1791--1803, 2006.

\bibitem[But03]{Butcher}
J.C. Butcher.
\newblock {\em Numerical {M}ethods for {O}rdinary {D}ifferential {E}quations}.
\newblock John Wiley and Sons Ltd., 2003.

\bibitem[CCGU12]{Gwiazda2}
J.A. Carrillo, R.M. Colombo, P.~Gwiazda, and A.~Ulikowska.
\newblock Structured populations, cell growth and measure valued balance laws.
\newblock {\em Journal of Differential Equations}, 252:3245--3277, 2012.

\bibitem[CCR11]{CanCarRos}
J.A. Ca{\~n}izo, J.A. Carrillo, and J.~Rosado.
\newblock A well-posedness theory in measures for some kinetic models of
  collective motion.
\newblock {\em Mathematical Models and Methods in Applied Sciences},
  21(3):515--539, 2011.

\bibitem[CCS15]{Colombo_Crippa_Spirito}
M.~Colombo, G.~Crippa, and S.~Spirito.
\newblock Renormalized solutions to the continuity equation with an integrable
  damping term.
\newblock {\em Calc. Var.}, 54(2):1831--1845, 2015.

\bibitem[CDF{\etalchar{+}}11]{CarrDiFranFigLauSlep11}
J.A. Carrillo, M.~DiFrancesco, A.~Figalli, T.~Laurent, and D.~Slep\v{c}ev.
\newblock Global-in-time weak measure solutions and finite-time aggregation for
  nonlocal interaction equations.
\newblock {\em Duke Math. J.}, 156:229--271, 2011.

\bibitem[CFRT10]{CarrilloFornasierRosadoToscani}
J.A. Carrillo, M.~Fornasier, J.~Rosado, and G.~Toscani.
\newblock Asymptotic flocking dynamics for the kinetic {C}ucker-{S}male model.
\newblock {\em SIAM J. Math. Anal.}, 42(1):218--236, 2010.

\bibitem[CG09]{ColomboGuerra}
R.M. Colombo and G.~Guerra.
\newblock Differential equations in metric spaces with applications.
\newblock {\em Discrete and Continuous Dynamical Systems A}, 23(3):733--753,
  2009.

\bibitem[CLM13]{CrippaLecureux}
G.~Crippa and M.~L\'{e}cureux-Mercier.
\newblock Existence and uniqueness of measure solutions for a system of
  continuity equations with non-local flow.
\newblock {\em Nonlinear Differential Equations and Applications NoDEA},
  20(3):523--537, 2013.

\bibitem[CPT14]{CrisPiccTosinBook}
E.~Cristiani, B.~Piccoli, and A.~Tosin.
\newblock {\em Multiscale {M}odeling of {P}edestrian {D}ynamics}, volume~12 of
  {\em Modeling, Simulation \& Applications}.
\newblock Springer International Publishing Switzerland, 2014.

\bibitem[DG05]{DG}
O.~Diekmann and Ph. Getto.
\newblock Boundedness, global existence and continuous dependence for nonlinear
  dynamical systems describing physiologically structured populations.
\newblock {\em Journal of Differential Equations}, 215(2):268--319, 2005.

\bibitem[DU77]{Diestel-Uhl}
J.~Diestel and J.J. {Uhl jr.}
\newblock {\em Vector {M}easures}.
\newblock Amer. Math. Soc., Providence, 1977.

\bibitem[Dud66]{Dudley1}
R.M. Dudley.
\newblock Convergence of {B}aire measures.
\newblock {\em Stud. Math.}, 27:251--268, 1966.

\bibitem[Dud74]{Dudley2}
R.M. Dudley.
\newblock Correction to: ``{C}onvergence of {B}aire measures".
\newblock {\em Stud. Math.}, 51:275, 1974.

\bibitem[Dud04]{Dudley}
R.M. Dudley.
\newblock {\em Real {A}nalysis and {P}robability}.
\newblock Cambridge University Press, 2004.

\bibitem[EHM15a]{EvHiMu_JDE}
J.H.M. Evers, S.C. Hille, and A.~Muntean.
\newblock Mild solutions to a measure-valued mass evolution problem with flux
  boundary conditions.
\newblock {\em Journal of Differential Equations}, 259:1068--1097, 2015.

\bibitem[EHM15b]{EversMBE}
J.H.M. Evers, S.C. Hille, and A.~Muntean.
\newblock Modelling with measures: {A}pproximation of a mass-emitting object by
  a point source.
\newblock {\em Mathematical Biosciences and Engineering}, 12(2):357--373, 2015.

\bibitem[Eve15]{Evers_PhD}
J.H.M. Evers.
\newblock {\em {E}volution {E}quations for {S}ystems {G}overned by {S}ocial
  {I}nteractions}.
\newblock PhD thesis, Eindhoven University of Technology, 2015.

\bibitem[FM53]{Fortet-Mourier:1953}
R.~Fortet and E.~Mourier.
\newblock Convergence de la r\'epartition empirique vers la r\'epartition
  th\'eorique.
\newblock {\em Ann. Sci. E.N.S.}, 70(3):276--285, 1953.

\bibitem[GHS{\etalchar{+}}14]{Schleper}
S.~G\"{o}ttlich, S.~Hoher, P.~Schindler, V.~Schleper, and A.~Verl.
\newblock Modeling, simulation and validation of material flow on conveyor
  belts.
\newblock {\em Appl. Math. Modell.}, 38:3295--3313, 2014.

\bibitem[GJMC12]{Gwiazda-Jamroz_ea:2012}
P.~Gwiazda, G.~Jamr\'oz, and A.~Marciniak-Czochra.
\newblock Models of discrete and continuous cell differentiation in the
  framework of transport equation.
\newblock {\em SIAM J. Math. Anal.}, 44(2):1103--1133, 2012.

\bibitem[GLMC10]{Gwiazda1}
P.~Gwiazda, T.~Lorenz, and A.~Marciniak-Czochra.
\newblock A nonlinear structured population model: Lipschitz continuity of
  measure-valued solutions with respect to model ingredients.
\newblock {\em Journal of Differential Equations}, 248:2703--2735, 2010.

\bibitem[Hoo13]{Hoogwater-thesis}
R.~Hoogwater.
\newblock {N}on-linear {S}tructured {P}opulation {M}odels: {A}n {A}pproach with
  {S}emigroups on {M}easures and {E}uler's {M}ethod.
\newblock Master's thesis, Leiden University, February 2013.

\bibitem[HW09]{Hille-Worm:SF}
S.C. Hille and D.T.H. Worm.
\newblock Continuity properties of {M}arkov semigroups and their restrictions
  to invariant {$L^1$}-spaces.
\newblock {\em Semigroup Forum}, 79:575--600, 2009.

\bibitem[LMS02]{Lasota-Myjak-Szarek}
A.~Lasota, J.~Myjak, and T.~Szarek.
\newblock Markov operators with a unique invariant measure.
\newblock {\em J. Math. Anal. Appl.}, 276:343--356, 2002.

\bibitem[PR13]{PiccoliRossi}
B.~Piccoli and F.~Rossi.
\newblock Transport equation with nonlocal velocity in {W}asserstein spaces:
  convergence of numerical schemes.
\newblock {\em Acta Applicandae Mathematicae}, 124(1):73--105, 2013.

\bibitem[TF11]{TosinFrasca}
A.~Tosin and P.~Frasca.
\newblock Existence and approximation of probability measure solutions to
  models of collective behaviors.
\newblock {\em Networks and Heterogeneous Media}, 6(3):561--596, 2011.

\bibitem[vMM14]{vMeurs14}
P.~van Meurs and A.~Muntean.
\newblock Upscaling of the dynamics of dislocation walls.
\newblock {\em Advances in Mathematical Sciences and Applications},
  24(2):401--414, 2014.

\bibitem[\v{S}94]{Sikic}
H.~\v{S}iki\'c.
\newblock Nonlinear perturbations of positive semigroups.
\newblock {\em Semigroup Forum}, 48:273--302, 1994.

\bibitem[Zah00]{Zaharopol}
R.~Zaharopol.
\newblock Fortet-{M}ourier norms associated with some iterated function
  systems.
\newblock {\em Stat. Prob. Letters}, 50:149--154, 2000.

\end{thebibliography}

\end{document}